\documentclass{amsart}

\usepackage{amsmath}%
\usepackage{amsfonts}%
\usepackage{amssymb}%
\usepackage{graphicx}

\newtheorem{theo}{Theorem}[section]
\newtheorem{lemme}[theo]{Lemma}

\newtheorem{defi}[theo]{Definition}

\newtheorem{prop}[theo]{Proposition}
\newtheorem{rmq}[theo]{Remark}

\usepackage{color}
\definecolor{violet}{rgb}{0.46,0.39,0.87}

\numberwithin{equation}{section}
\usepackage{marvosym}

\DeclareMathAlphabet{\mathonebb}{U}{bbold}{m}{n}
\newcommand{\one}{\ensuremath{\mathonebb{1}}}
\begin{document}
\title[Sub-critical Keller-Segel]{Propagation of chaos for a sub-critical Keller-Segel model}
\author{David Godinho, Cristobal Quininao}

\begin{abstract}
This paper deals with a sub-critical Keller-Segel equation. Starting from the stochastic particle system associated with it, we show well-posedness results and the propagation of chaos property. More precisely, we show that the empirical measure of the system tends towards the unique solution of the limit equation as the number of particles goes to infinity.
\end{abstract}
\maketitle



David Godinho \Letter: Laboratoire d'Analyse et de Math\'ematiques Appliqu\'ees,
CNRS UMR 8050, Universit\'e Paris-Est, 61 avenue du G\'en\'eral de Gaulle, 94010 Cr\'eteil Cedex, France.

\textit{E-mail address:} david.godinho-pereira@u-pec.fr

Cristobal Qui\~ninao \Letter: The Mathematical Neuroscience Laboratory, CIRB and INRIA Bang Laboratory, Coll\`ege de France, 11, place Marcelin Berthelot, 75005 Paris, France.

\textit{E-mail address:} cristobal.quininao@college-de-france.fr

Cristobal Qui\~ninao: UPMC Univ Paris 06, UMR 7598, Laboratoire Jacques-Louis Lions, F-75005, Paris, France.
\vspace{1cm}

\section{Introduction and main results}

The subject of this paper is the convergence of a stochastic particle system to a non linear and non local equation which can be seen as a sub-critical version of the classical Keller-Segel equation.

\subsection{The sub-critical Keller-Segel Equation}  Consider the equation:
\begin{equation}\label{eq:KS}
 \frac{\partial f_t(x)}{\partial t}=\chi\,\nabla_x\cdot((K\ast f_t)(x))f_t(x))+\triangle_x f_t(x),
\end{equation}
where $f:\mathbb R_+\times\mathbb R^2\rightarrow\mathbb R$ and $\chi>0$. The force field kernel $K:\mathbb R^2\rightarrow\mathbb R^2$ comes from an attractive potential $\Phi:\mathbb R^2\rightarrow\mathbb R$ and is defined by
\begin{align} \label{KKellerSegel}
K(x):=\frac{x}{|x|^{\alpha+1}}=-\nabla\underbrace{\left(\frac1{\alpha-1}|x|^{1-\alpha}\right)}_{\Phi(x)}, \quad \alpha\in(0,1).
\end{align}

The standard Keller-Segel equation correspond to the critical case $K(x)=x/|x|^2$ (i.e., more singular at $x=0$) and it describes a model of chemotaxis, i.e., the movement of cells (usually bacteria or amoebae) which are attracted by some chemical substance that they product. This equation has been first introduced by Keller and Segel in \cite{KS1,KS2}. Blanchet-Dolbeault-Perthame showed in \cite{BDP} some nice results on existence of global weak solutions if the parameter $\chi$ (which is the sensitivity of the bacteria to the chemo-attractant) is smaller than $8\pi/M$ where $M$ is the initial mass (here $M$ will always be 1 since we will deal with probability measures). For more details on the subject, see \cite{HOR1,HOR2}.

\subsection{The particle system}
We consider the following system of particles
\begin{align} \label{systpart}
\forall i=1,...,N, \quad X_t^{i,N}=X_0^{i,N}-\frac{\chi}{N}\sum_{j=1,j\neq i}^N\int_0^tK(X_s^{i,N}-X_s^{j,N})ds+\sqrt2 B_t^i,
\end{align}
where $(B^i)_{i=1,...,N}$ is an independent family of 2D standard Brownian motions and $K$ is defined in \eqref{KKellerSegel}. We will show in the sequel that there is propagation of chaos to the solution of the following nonlinear S.D.E linked with \eqref{eq:KS} (see the next paragraph)
\begin{align} \label{edsKS}
X_t=X_0-\chi\int_0^t\int_{\mathbb R^2}K(X_s-x)f_s(dx)ds+\sqrt2 B_t,
\end{align}
where $f_t=\mathcal{L}(X_t)$.

\subsection{Weak solution for the P.D.E}
For any Polish space $E$, we denote by $\textbf{P}(E)$ the set of all probability measures on $E$ which we endow with the topology of weak convergence defined by duality against functions of $C_b(E)$. We give the notion of weak solution that we use in this paper.

\begin{defi}\label{def:Weak}
 We say that $f=(f_t)_{t\geq0}\in C([0,\infty),{\bf P} (\mathbb R^2))$ is a weak solution to \eqref{eq:KS} if
 \begin{equation}\label{eq:L1bound}
  \forall\, T>0,\qquad\int_0^T\int_{\mathbb R^2}\int_{\mathbb R^2}|K(x-y)|\,f_s(dx)\,f_s(dy)\,ds<\infty,
 \end{equation}
 and if for all $\varphi\in C_b^2(\mathbb R^2)$, all $t\geq0$,
 \begin{align}\label{eq:weak}
 \nonumber \int_{\mathbb R^2}\varphi(x)f_t(dx)=&\int_{\mathbb R^2}\varphi(x)f_0(dx)+\int_0^t\int_{\mathbb R^2}\triangle_x\varphi(x)f_s(dx)\,ds\\
 &-\chi\,\int_0^t\int_{\mathbb R^2}\int_{\mathbb R^2}K(x-y)\cdot\nabla_x\varphi(x)f_s(dy)f_s(dx)\,ds.
 \end{align}
\end{defi}

\begin{rmq} \label{rmq:weaksol}
 We can see easily that if $(X_t)_{t\geq0}$ is a solution to \eqref{edsKS}, then setting $f_t=\mathcal L(X_t)$ for any $t\geq0$, $(f_t)_{t\geq0}$ is a weak solution of \eqref{eq:KS} in the sense of Definition \ref{def:Weak} provided it satisfies \eqref{eq:L1bound}. Indeed, by It\^o's formula, we find that for $\varphi\in C_b^2(\mathbb R^2)$, 
 \begin{align}
  \nonumber \varphi(X_t)=&\varphi(X_0)-\chi\int_0^t\nabla_x\varphi(X_s)\cdot\int_{\mathbb R^2}K(X_s-y)f_s(dy)\,ds\\
 \nonumber &+\int_0^t\sqrt{2}\nabla_x\varphi(X_s)\cdot dB_s+\int_0^t\triangle_x\varphi(X_s)ds.
 \end{align}
 Taking expectations, we get \eqref{eq:weak}.
\end{rmq}

\subsection{Notation and propagation of chaos} For $N\geq2$, we denote by $\textbf{P}_{sym}(E^N)$ the set of symmetric probability measures on $E^N$, i.e. the set of probability measures which are laws of exchangeable $E^N$-valued random variables. 

We consider for any $F\in\textbf{P}_{sym}((\mathbb R^2)^N)$ with a density (a finite moment of positive order is also required in order to define the entropy) the Boltzmann entropy and the Fisher information which are defined by
\begin{align*}
H(F):=\frac{1}{N}\int_{(\mathbb R^2)^N}F(x)\log F(x)dx\quad {\rm and}\quad I(F):=\frac{1}{N}\int_{(\mathbb R^2)^N}\frac{|\nabla F(x)|^2}{F(x)}dx.
\end{align*} 
We also define ($x_i\in\mathbb R^2$ stands for the i-th coordinate of $x\in(\mathbb R^2)^N$), for $k\geq0$, 
\begin{align*}
M_k(F):=\frac{1}{N}\int_{(\mathbb R^2)^N}\sum_{i=1}^N|x_i|^kF(dx).
\end{align*}
Observe that we proceed to the normalization by $1/N$ in order to have, for any $f\in\textbf{P}(\mathbb R^2)$,
\begin{align*}
H(f^{\otimes N})=H(f),\quad I(f^{\otimes N})=I(f) \quad {\rm and}\quad M_k(f^{\otimes N})=M_k(f).
\end{align*}
We introduce the space $\textbf{P}_1(\mathbb R^2):=\{f\in\textbf{P}(\mathbb R^2)$, $M_1(f)<\infty$\} and we recall the definition of the Wasserstein distance: if $f,g\in\textbf{P}_1(\mathbb R^2)$,
\begin{align*}
\mathcal W_1(f,g)=\inf\Big\{\int_{\mathbb R^2\times\mathbb R^2}|x-y|R(dx,dy)\Big\},
\end{align*}
where the infimum is taken over all probability measures $R$ on $\mathbb R^2\times\mathbb R^2$ with $f$ for first marginal and $g$ for second marginal. It is known that the infimum is reached. See e.g. Villani \cite{VIL2} for many details on the subject.

We now define the notion of propagation of chaos.

\begin{defi}
Let $X$ be some $E$-valued random variable. A sequence $(X_1^N,...,X_N^N)$ of exchangeable $E$-valued random variables is said to be $X$-chaotic if one of the three following equivalent conditions is satisfied:\\
(i) $(X_1^N,X_2^N)$ goes in law to 2 independent copies of $X$ as $N\rightarrow+\infty$;\\
(ii) for all $j\geq1$, $(X_1^N,...,X_j^N)$ goes in law to $j$ independent copies of $X$ as $N\rightarrow+\infty$;\\
(iii) the empirical measure $\mu_{X^N}^N:=\frac{1}{N}\sum_{i=1}^N \delta_{X_i^N}\in\textbf{P}(E)$ goes in law to the constant $\mathcal L(X)$ as $N\rightarrow+\infty$.
\end{defi}
We refer to \cite{SZN} for the equivalence of the three conditions or \cite[Theorem 1.2]{HM} where the equivalence is established in a quantitative way.

Propagation of chaos in the sense of Sznitman holds for a system of $N$ exchangeable particles evolving in time if when the initial conditions $(X_0^{1,N},...,X_0^{N,N})$ are $X_0$-chaotic, the trajectories $((X_t^{1,N})_{t\geq0},...,(X_t^{N,N})_{t\geq0})$ are $(X_t)_{t\geq0}$-chaotic, where $(X_t)_{t\geq0}$ is the (unique) solution of the expected (one-particle) limit model.  

We finally recall a stronger (see \cite{HM}) sense of chaos introduced by Kac in \cite{KAC} and formalized recently in \cite{CCLLV}: the entropic chaos.
\begin{defi}
Let $f$ be some probability measure on $E$. A sequence $(F^N)$ of symmetric probability measures on $E^N$ is said to be entropically $f$-chaotic if 
\begin{align*}
F_1^N\rightarrow f \quad{\rm weakly\  in\  }{\bf P}(E)\quad{\rm and} \quad H(F^N)\rightarrow H(f)\quad {\rm as}\quad N\rightarrow\infty,
\end{align*}
where $F_1^N$ stands for the first marginal of $F^N$.
\end{defi}
We can observe that since the entropy is lower semi continuous (so that $H(f)\leq\liminf_N H(F^N)$) and is convex, the entropic chaos (which requires $\lim_N H(F^N)=H(f)$) is a stronger notion of convergence which implies that for all $j\geq1$, the density of the law of $(X_1^N,...,X_j^N)$ goes to $f^{\otimes j}$ strongly in $L^1$ as $N\rightarrow\infty$ (see \cite{BRE}). 

\subsection{Main results}

We first give a result of existence and uniqueness for \eqref{eq:KS}.
\begin{theo}\label{thm:ExUnEDP}
Let $\alpha\in(0,1)$. Assume that $f_0\in \textbf{P}_1(\mathbb R^2)$ is such that $H(f_0)<\infty$. \\
 (i) There exists a unique weak solution $f$ to \eqref{eq:KS} such that
 \begin{equation}\label{eq:FinFirstMom}
  f\in L_{loc}^\infty([0,\infty),\textbf{P}_1(\mathbb R^2))\cap L_{loc}^1([0,\infty);L^p(\mathbb R^2)) \quad {\rm for\ some\ } p>\frac{2}{1-\alpha}.
 \end{equation}
   (ii) This solution furthermore satisfies that for all $T>0$,
 \begin{equation}\label{eq:Fishinff}
  \int_0^TI(f_s)ds<\infty,
 \end{equation}
 for any $q\in[1,2)$ and for all $T>0$,
 \begin{equation}\label{eq:Derivative}
   \nabla_x f\in L^{2q/(3q-2)}(0,T;L^q(\mathbb R^2)),
  \end{equation}
  for any $p\geq1$,
  \begin{equation}\label{eq:Regularity}
  f\in C([0,\infty);L^1(\mathbb R^2))\cap C((0,\infty);L^p(\mathbb R^2)),
 \end{equation}
 and that for any $\beta\in C^1(\mathbb R)\cap W^{2,\infty}_{loc}(\mathbb R)$ such that $\beta^{\prime\prime}$ is piecewise continuous and vanishes outside a compact set,
 \begin{align}\label{eq:Renormalization}
   \partial_t\beta(f)=&\chi\,(K\ast f)\cdot\nabla_x(\beta(f))+\triangle_x\beta(f)\\
  \nonumber&-\beta^{\prime\prime}(f)|\nabla_x f|^2+\chi\,\beta^\prime(f_s)f_s(\nabla_x\cdot K\ast f_s),
 \end{align}
 on $[0,\infty)\times\mathbb R^2$ in the distributional sense.
\end{theo}

We denote by $F_0^N$ the law of $(X_0^{i,N})_{i=1,...,N}$. We assume that for some $f_0\in \textbf{P}(\mathbb{R}^2)$,
\begin{align} \label{f0chaotic}
\begin{cases}
F_0^N\in\textbf{P}_{sym}((\mathbb{R}^2)^N) \quad {\rm is }\quad f_0-{\rm chaotic};\\
\displaystyle \sup_{N\geq2}M_1(F_0^N)<\infty,\quad \sup_{N\geq2}H(F_0^N)<\infty.
\end{cases}
\end{align}
Observe that this condition is satisfied if the random variables $(X_0^{i,N})_{i=1,...,N}$ are i.i.d. with law $f_0\in \textbf{P}_1(\mathbb R^2)$ such that $H(f_0)<\infty$. The next result states the well-posedness for the particle system \eqref{systpart}.

\begin{theo} \label{HMIFt}
Let $\alpha\in(0,1)$.\\
(i) Let $N\geq2$ be fixed and assume that $M_1(F_0^N)<\infty$ and $H(F_0^N)<\infty$. There exists a unique strong solution $(X_t^{i,N})_{t\geq0,i=1,...,N}$ to (\ref{systpart}). Furthermore, the particles a.s. never collapse i.e. it holds that a.s., for any $t\geq0$ and $i\neq j$,  $X_t^{i,N}\neq X_t^{j,N}$.\\
(ii) Assume \eqref{f0chaotic}. If for all $t\geq0$, we denote by $F_t^{N}\in \textbf{P}_{sym}((\mathbb R^2)^N)$ the law of $(X_t^{i,N})_{i=1,...,N}$, then there exist a constant C depending on $\chi$, $\sup_{N\geq2}H(F_0^N)$ and $\sup_{N\geq2}M_1(F_0^N)$ such that for all $t\geq0$ and $N\geq2$
\begin{align*}
H(F_t^{N})\leq C(1+t), \quad M_1(F_t^{N})\leq C(1+t), \quad \int_0^t I(F_s^{N})ds\leq C(1+t).
\end{align*}
Furthermore for any $T>0$,
\begin{align}\label{EsupX1}
\mathbb E\Big[\sup_{t\in[0,T]}|X_t^{1,N}|\Big]\leq C(1+T).
\end{align}
We also have
\begin{align} \label{entroFishFt}
H(F_t^N)+\int_0^tI(F_s^N)ds\leq H(F_0^N)+\frac{\chi}{N^2}\sum_{i\neq j}\int_0^t\mathbb E\big[{\rm div} K(X_s^{i,N}-X_s^{j,N})\big]\,ds.
\end{align}
\end{theo}

We next state a well-posedness result for the nonlinear S.D.E. \eqref{edsKS}.

\begin{theo} \label{exuniedsKS}
Let $\alpha\in(0,1)$ and $f_0\in \textbf{P}_1(\mathbb{R}^2)$ such that $H(f_0)<\infty$. There exists a unique strong solution $(X_t)_{t\geq0}$ to (\ref{edsKS}) such that for some $p>2/(1-\alpha)$,
\begin{align} \label{Fishmom}
(f_t)_{t\geq0}\in L_{loc}^\infty([0,\infty),\textbf{P}_1(\mathbb R^2))\cap L_{loc}^1([0,\infty);L^p(\mathbb R^2)),
\end{align}
where $f_t$ is the law of $X_t$. Furthermore, $(f_t)_{t \geq0}$ is the unique solution to \eqref{eq:KS} given in Theorem \ref{thm:ExUnEDP}.
\end{theo}

We finally give the result about propagation of chaos.

\begin{theo} \label{propchaos}
Let $\alpha\in(0,1)$. Assume (\ref{f0chaotic}). For each $N\geq2$, consider the unique solution $(X_t^{i,N})_{i=1,...,N,t\geq0}$ to (\ref{systpart}). Let $(X_t)_{t\geq0}$ be the unique solution to (\ref{edsKS}).\\
(i) The sequence $(X_t^{i,N})_{i=1,...,N,t\geq0}$ is $(X_t)_{t\geq0}$-chaotic. In particular, the empirical measure $Q^N:=\frac{1}{N}\sum_{i=1}^N\delta_{(X_t^{i,N})_{t\geq0}}$ goes in law to $\mathcal L((X_t)_{t\geq0})$ in $\textbf{P}(C((0,\infty),\mathbb R^2))$.  \\
(ii) Assume furthermore that $\lim_N H(F_0^N)=H(f_0)$. For all $t\geq0$, the sequence $(X_t^{i,N})_{i=1,...,N}$ is then $X_t$-entropically chaotic. In particular, for any $j\geq1$ and any $t\geq0$, denoting by $F_{tj}^N$ the density of the law of $(X_t^{1,N},...,X_t^{j,N})$, it holds that
\begin{align*}
\lim_{N\rightarrow\infty} ||F_{tj}^N-f_t^{\otimes j}||_{L^1((\mathbb R^2)^j)}=0.
\end{align*}
\end{theo}
We can observe that the condition $\lim_N H(F_0^N)=H(f_0)$ is satisfied if the random variables $(X_0^{i,N})_{i=1,...,N}$ are i.i.d. with law $f_0$ such that $H(f_0)<\infty$.
\subsection{Comments}
This paper is some kind of adaptation of the work of Fournier-Hauray-Mischler in \cite{FHM} where they show the propagation of chaos of some particle system for the 2D viscous vortex model. We use the same methods for a sub-critical Keller-Segel equation. The proofs are thus sometimes very similar to those in  \cite{FHM} but there are some differences due to the facts that i) there are no circulation parameter ($\mathcal M_i^N$ in \cite{FHM}): this simplify the situation since we thus deal with solutions which are probabilities and ii) the kernel is not the same: it is not divergence-free and we thus have to deal with some additional terms in our computations. We can also notice that due to this fact, we have no already known result for the existence and uniqueness of the particle system that we consider. 

The proof of Theorem \ref{thm:ExUnEDP} follows the ideas of renormalisation solutions to a PDE introduced by Di Perna and Lions in \cite{dPL} and developed since then. The key point is to be able to find good \textit{a priori} estimates which allow us to approximate the weak solutions by regular functions, i.e., to use $C^k$ functions instead of $L^1$. Then, using these estimates, one can pass to the limit and go back to the initial problem. One can further see that the uniqueness result is proved based on coupling methods and the Wasserstein distance. This will allow us to use more general initial conditions than we could use in a strictly deterministic framework.

The proof of existence and uniqueness for the particle system \eqref{systpart} (Theorem \ref{HMIFt}) use some nice arguments. Like for S.D.Es with locally Lipschitz coefficients, we show existence and uniqueness up to an explosion time and the interesting part of the proof is to show that this explosion time is infinite a.s.

To our knowledge, there is no other work that give a convergence result of some particle system for a chemotaxis model with a singular kernel $K$ and without cutoff parameter. In \cite{Stevens2}, Stevens studies a particle system with two kinds of particles corresponding to bacteria and chemical substance. She shows convergence of the system for smooth initial data (lying in $C_b^3(\mathbb R^d)$) and for regular kernels (continuously differentiable and bounded together with their derivatives). In \cite{Schmeiser2}, Haskovec and Schmeiser consider a kernel with a cutoff parameter $K_\epsilon(x)=\frac{x}{|x|(|x|+\epsilon)}$. They get some well-posedness result for the particle system and they show the weak convergence of subsequences due to a tightness result (observe that here we have propagation of chaos and also entropic chaos). In a recent work \cite{Calvez13}, Calvez and Corrias work on some one-dimensional Keller-Segel model. They study a dynamical particle system  for which they give a global existence result under some assumptions on the initial distribution of the particles that prevents collisions. They also give two blow-up criteria for the particle system they do not state a convergence  result for this system.    

Finally, it is important to notice that the present method can not be directly adapted for the standard case $\alpha=1$ because in this last situation the entropy and the Fisher information are not controlled. 

\subsection{Plan of the paper}
In the next section, we give some preliminary results. In Section 3, we establish the well-posedness of the particle system \eqref{systpart}. In Section~4, we prove the tightness of the particle system and we show that any limit point belongs to the set of solutions to the nonlinear S.D.E. \eqref{edsKS}. In Section~5, we show that the P.D.E. \eqref{eq:KS} and the nonlinear S.D.E. \eqref{edsKS} are well-posed and we show the propagation of chaos. Finally, in the last section, we improve the regularity of the solution, give some renormalization results for the solution to \eqref{eq:KS} and we conclude with the entropic chaos.

\section{Preliminaries}

In this section, we recall some lemmas stated in \cite{FHM} and \cite{HM} and we state a result on the regularity of the kernel $K$ defined in \eqref{KKellerSegel}. The first result tells us that pairs of particles which law have finite Fisher information cannot be too close.

\begin{lemme} \label{lemFish} $($\cite[Lemma 3.3]{FHM}$)$
Consider $F\in \textbf{P}(\mathbb R^2\times\mathbb R^2)$ with finite Fisher information and $(X_1,X_2)$ a random variable with law $F$. Then for any $\gamma\in(0,2)$ and any $\beta>\gamma/2$ there exists $C_{\gamma,\beta}$ so that
\begin{align*}
\mathbb E[|X_1-X_2|^{-\gamma}]=\int_{\mathbb R^2\times\mathbb R^2}\frac{F(x_1,x_2)}{|x_1-x_2|^\gamma}dx_1dx_2\leq C_{\gamma,\beta}(I(F)^\beta+1).
\end{align*}
\end{lemme}
In the next lemma, we see that the Fisher information of the marginals of some $F\in \textbf{P}_{sym}((\mathbb R^2)^N)$ is smaller than the Fisher information of $F$.

\begin{lemme} \label{lemFishmarg} $($\cite[Lemma 3.7]{HM}$)$
For any $F\in \textbf{P}_{sym}((\mathbb R^2)^N)$ and $1\leq l\leq N$, $I(F_l)\leq I(F)$, where $F_l\in \textbf{P}_{sym}((\mathbb R^2)^l)$ denotes the marginal probability of $F$ on the $l$-th block of variables.
\end{lemme}
The following lemma allows us to control from below the entropy of some $F\in \textbf{P}_k((\mathbb R^2)^N)$ by its moment of order $k$ for any $k>0$. 

\begin{lemme} \label{lementromoment} $($\cite[Lemma 3.1]{FHM}$)$
For any $k,\lambda\in(0,\infty)$, there is a constant $C_{k,\lambda}\in\mathbb R$ such that for any $N\geq1$, any $F\in \textbf{P}_k((\mathbb R^2)^N)$,
\begin{align*}
H(F)\geq -C_{k,\lambda}-\lambda M_k(F).
\end{align*}
\end{lemme}

The next result tells us that a probability measure on $\mathbb R^2$ with finite Fisher information belongs to $L^p$ for any $p\geq1$ and its derivatives, to $L^q$ for any $q\in[1.2)$.
\begin{lemme} \label{lemLpFish} $($\cite[Lemma 3.2]{FHM}$)$
For any $f\in \textbf{P}(\mathbb R^2)$ with finite Fisher information, there holds
\begin{align*}
&\forall p\in[1,\infty),\quad \|f\|_{L^p(\mathbb R^2)}\leq C_p I(f)^{1-1/p},\\
&\forall q\in[1,2),\quad \|\nabla_xf\|_{L^q(\mathbb R^2)}\leq C_q I(f)^{3/2-1/q}.
\end{align*}
\end{lemme}

We end this section with the following result on $K$.
\begin{lemme} \label{Kx-Ky}
Let $\alpha\in(0,1)$. There exists a constant $C_\alpha$ such that for all $x,y\in \mathbb R^2$
\begin{align*}
|K(x)-K(y)|\leq C_\alpha |x-y|\Big(\frac{1}{|x|^{\alpha+1}}+\frac{1}{|y|^{\alpha+1}}\Big).
\end{align*}
\end{lemme}

\textbf{Proof.}
We have
\begin{align*}
|K(x)-K(y)|&=\Big|x\Big(\frac{1}{|x|^{\alpha+1}}-\frac{1}{|y|^{\alpha+1}}\Big)+\frac{x-y}{|y|^{\alpha+1}}\Big|\\
&\leq |x||x-y|(\alpha+1)\max\Big(\frac{1}{|x|^{\alpha+2}},\frac{1}{|y|^{\alpha+2}}\Big)+\frac{|x-y|}{|y|^{\alpha+1}}.
\end{align*}
By symmetry, we also have 
\begin{align*}
|K(x)-K(y)|&\leq |y||x-y|(\alpha+1)\max\Big(\frac{1}{|x|^{\alpha+2}},\frac{1}{|y|^{\alpha+2}}\Big)+\frac{|x-y|}{|x|^{\alpha+1}}.
\end{align*}
So we deduce that
\begin{align*}
|K(x)-K(y)|&\leq |x-y|\Big[(\alpha+1)\min(|x|,|y|)\max\Big(\frac{1}{|x|^{\alpha+2}},\frac{1}{|y|^{\alpha+2}}\Big)\\
&\quad+\frac{1}{|x|^{\alpha+1}}+\frac{1}{|y|^{\alpha+1}}\Big]\\
&\leq |x-y|\Big[
(\alpha+1)\frac{1}{\min(|x|,|y|)^{\alpha+1}}+\frac{1}{|x|^{\alpha+1}}+\frac{1}{|y|^{\alpha+1}}\Big]\\
&\leq 2(\alpha+2)|x-y|\Big(\frac{1}{|x|^{\alpha+1}}+\frac{1}{|y|^{\alpha+1}}\Big).
\end{align*}
which concludes the proof. \hfill$\square$

\section{Well-posedness for the system of particles} \label{sect:WPsystpart}

Let's now introduce another particle system with a regularized kernel. We set, for $\epsilon\in(0,1)$,
\begin{align}  \label{Keps}
K_\epsilon(x)=\frac{x}{\max(|x|,\epsilon)^{\alpha+1}},
\end{align}
which obviously satisfies $|K_\epsilon(x)-K_\epsilon(y)|\leq C_{\alpha,\epsilon} |x-y|$ and we consider the following system of S.D.E.s 
\begin{align} \label{systparteps}
\forall i=1,...,N, \quad X_t^{i,N,\epsilon}=X_0^{i,N}-\frac{\chi}{N}\sum_{j=1,j\neq i}^N\int_0^tK_\epsilon(X_s^{i,N,\epsilon}-X_s^{j,N,\epsilon})ds+\sqrt2 B_t^i,
\end{align}
for which strong existence and uniqueness thus holds.

The following result will be useful for the proof of Theorem \ref{HMIFt}. Its proof is very similar to the proof of \cite[Proposition 5.1]{FHM}.
\begin{prop} \label{HMIFteps}
Let $\alpha\in(0,1)$.\\
(i) Let $N\geq2$ be fixed. Assume that $M_1(F_0^N)<\infty$ and $H(F_0^N)<\infty$. For all $t\geq0$, we denote by $F_t^{N,\epsilon}\in \textbf{P}_{sym}((\mathbb R^2)^N)$ the law of $(X_t^{i,N,\epsilon})_{i=1,...,N}$. Then 
\begin{align} \label{entroFteps}
H(F_t^{N,\epsilon})=&H(F_0^N)+\frac{\chi}{N^2}\sum_{i\neq j}\int_0^t\int_{(\mathbb R^2)^N} {\rm div} K_\epsilon(x_i-x_j)F_s^{N,\epsilon}(x)dsdx\\
\nonumber
&-\int_0^tI(F_s^{N,\epsilon})ds.
\end{align}
(ii) There exists a constant C which depends on $\chi$, $H(F_0^N)$ and $M_1(F_0^N)$ such that for all $t\geq0$ and $N\geq2$,
\begin{align} \label{eqHMIFteps}
H(F_t^{N,\epsilon})\leq C(1+t), \quad M_1(F_t^{N,\epsilon})\leq C(1+t), \quad \int_0^t I(F_s^{N,\epsilon})ds\leq C(1+t).
\end{align}
Furthermore,
\begin{align} \label{EsupXi2}
\mathbb E\Big[\sup_{[0,T]}|X_t^{1,N,\epsilon}|\Big]\leq C(1+T).
\end{align}
\end{prop}

\textbf{Proof.}
Let $\varphi\in C_b^2((\mathbb R^2)^N)$, and $t\geq0$ be fixed. Using It\^o's formula, we compute the expectation of $\varphi(X_t^{1,N,\epsilon},...,X_t^{N,N,\epsilon})$ and get (recall that $x_i\in\mathbb R^2$ stands for the i-th coordinate of $x\in(\mathbb R^2)^N$)
\begin{align} \label{Ftphi}
\frac{d}{dt}\int_{(\mathbb R^2)^N}\varphi(x)F_t^{N,\epsilon}(dx)=&-\frac{\chi}{N}\int_{(\mathbb R^2)^N}\sum_{i\neq j}K_\epsilon(x_i-x_j)\cdot\nabla_{x_i}\varphi(x)F_t^{N,\epsilon}(dx)\\
\nonumber&+\int_{(\mathbb R^2)^N}\triangle_x\varphi(x)F_t^{N,\epsilon}(dx).
\end{align}
We deduce that $F^{N,\epsilon}$ is a weak solution to 
\begin{align}\label{solfaibleFt}
\partial_t F_t^{N,\epsilon}(x)=\frac{\chi}{N}\sum_{i\neq j} {\rm div}_{x_i}(F_t^
{N,\epsilon}(x)K_\epsilon(x_i-x_j))+\triangle_x F_t^{N,\epsilon}(x).
\end{align}
We are now able to compute the evolution of the entropy.
\begin{align*}
\frac{d}{dt}H(F_t^{N,\epsilon})&=\frac{1}{N}\int_{(\mathbb R^2)^N}\partial_t F_t^{N,\epsilon}(x)(1+\log F_t^{N,\epsilon}(x))dx\\
&=\frac{\chi}{N^2}\sum_{i\neq j}\int_{(\mathbb R^2)^N} {\rm div}_{x_i}(F_t^{N,\epsilon}(x)K_\epsilon(x_i-x_j))(1+\log F_t^{N,\epsilon}(x))dx\\
&\quad+\frac{1}{N}\int_{(\mathbb R^2)^N}\triangle_x F_t^{N,\epsilon}(x)(1+\log F_t^{N,\epsilon}(x))dx.
\end{align*}
Performing some integrations by parts, we get
\begin{align*}
\frac{d}{dt}H(F_t^{N,\epsilon})&=-\frac{\chi}{N^2}\sum_{i\neq j}\int_{(\mathbb R^2)^N}K_\epsilon(x_i-x_j)\cdot\nabla_{x_i}F_t^{N,\epsilon}(x)dx-I(F_t^{N,\epsilon})\\
&=\frac{\chi}{N^2}\sum_{i\neq j}\int_{(\mathbb R^2)^N}{\rm div}K_\epsilon(x_i-x_j)F_t^{N,\epsilon}(x)dx-I(F_t
^{N,\epsilon}),
\end{align*}
and \eqref{entroFteps} follows. Using that div $K_\epsilon(x)=\frac{1-\alpha}{|x|^{\alpha+1}}\one_{\{|x|\geq\epsilon\}}+\frac{2}{\epsilon^{\alpha+1}}\one_{\{|x|<\epsilon\}}\leq \frac{2}{|x|^{\alpha+1}}$ and the exchangeability of the particles, we get
\begin{align*}
\frac{d}{dt}H(F_t^{N,\epsilon})&\leq\frac{2\chi}{N^2}\sum_{i\neq j}\int_{(\mathbb R^2)^N}\frac{F_t^{N,\epsilon}(x)}{|x_i-x_j|^{\alpha+1}}dx-I(F_t^{N,\epsilon})\\
\nonumber
&\leq 2\chi \int_{(\mathbb R^2)^N}\frac{F_t^{N,\epsilon}(x)}{|x_1-x_2|^{\alpha+1}}dx-I(F_t^{N,\epsilon}).
\end{align*}
Since $\alpha\in(0,1)$, we can use Lemma \ref{lemFish} with $\gamma=\alpha+1$ and $\beta$ such that $\frac{\alpha+1}{2}<\beta<1$, which gives
\begin{align*}
\int_{(\mathbb R^2)^N}\frac{F_t^{N,\epsilon}(x)dx}{|x_1-x_2|^{\alpha+1}}\leq C(I(F_{t2}^{N,\epsilon})^\beta+1),
\end{align*}
where
 $F_{t2}^{N,\epsilon}$ is the two-marginal of $F_t^{N,\epsilon}$. By Lemma \ref{lemFishmarg}, we have $I(F_{t2}^{N,\epsilon})\leq I(F_t^{N,\epsilon})$. Using that $Cx^\beta\leq C+\frac{x}{6\chi}$ (changing the value of the constant $C$), we thus get
\begin{align*}
\frac{d}{dt}H(F_t^{N,\epsilon})\leq C-\frac{2}{3}I(F_t^{N,\epsilon}),
\end{align*}
and thus
\begin{align} \label{entroFt}
H(F_t^{N,\epsilon})+\frac{2}{3}\int_0^t I(F_s^{N,\epsilon})ds\leq H(F_0^N)+Ct.
\end{align}
We now compute $M_1(F_t^{N,\epsilon})$. We first observe that 
\begin{align*}
M_1(F_t^{N,\epsilon})=\frac{1}{N}\int_{(\mathbb R^2)^N}\sum_{i=1}^N|x_i|F_t^{N,\epsilon}(dx)= \mathbb E[|X_t^{1,N,\epsilon}|],
\end{align*}
since the particles are exchangeable. We will need to control $\mathbb E[\sup_{[0,T]}|X_t^{1,N,\epsilon}|]$ in the sequel. We have
\begin{align} \label{EsupXi}
\mathbb E\Big[\sup_{[0,T]}|X_t^{1,N,\epsilon}|\Big]&\leq C\Big(\mathbb 
E[|X_0^1|]+\mathbb E\Big[\sup_{[0,T]}|B_t^1|\Big]\\
\nonumber
&\qquad+\mathbb E\Big[\sup_{t\in[0,T]}\Big|\frac{1}{N}\sum_{j\neq 1}\int_0^tK_\epsilon(X_s^{1,N,\epsilon}-X_s^{j,N,\epsilon})ds\Big|\Big]\Big)\\
\nonumber
&\leq C\Big(\mathbb E[|X_0^1|]+T+\frac{1}{N}\sum_{j\neq 1}\int_0^T \mathbb E[|K_\epsilon(X_s^{1,N,\epsilon}-X_s^{j,N,\epsilon})|]ds\Big)\\
\nonumber
&\leq C\Big(\mathbb E[|X_0^1|]+T+\int_0^T\mathbb E\Big[\frac{1}{|X_s^{1,N,\epsilon}-X_s^{2,N,\epsilon}|^\alpha}\Big]ds\Big).
\end{align} 
Using Lemma \ref{lemFish} with $\gamma=\alpha$ and $\beta$ such that $\frac{\alpha}{2}<\beta<1$ and recalling that $I(F_{t2}^{N,\epsilon})\leq I(F_t^{N,\epsilon})$, we get
\begin{align} \label{MkFt}
M_1(F_t^{N,\epsilon})&\leq C\Big(M_1(F_0^N)+T+\int_0^tI(F_t^{N,\epsilon})^\beta ds\Big)\\
\nonumber
&\leq C\Big(M_1(F_0^N)+T\Big)+\frac{1}{3}\int_0^tI(F_t^{N,\epsilon})ds,
\end{align}
where we used that $Cx^\beta\leq C+\frac{x}{3}$ (changing the value of $C$). Summing (\ref{entroFt}) and (\ref{MkFt}), we thus find
\begin{align*}
H(F_t^{N,\epsilon})+M_1(F_t^{N,\epsilon})+\frac{1}{3}\int_0^t I(F_s^{N,\epsilon})ds\leq H(F_0^N)+Ct+C(1+M_1(F_0^N)).
\end{align*}
Since the quantities $M_1$ and $I$ are positive, we immediately get $H(F_t^{N,\epsilon})\leq C(1+t)$. Using Lemma \ref{lementromoment}, we have $H(F_t^{N,\epsilon})\geq-C-M_1(F_t^{N,\epsilon})/2$, so that
\begin{align*}
M_1(F_t^{N,\epsilon})+\frac{1}{3}\int_0^t I(F_s^{N,\epsilon})ds\leq C(1+t)+M_1(F_t^{N,\epsilon})/2.
\end{align*}
Using again the positivity of $M_1$ and $I$, we easily get \eqref{eqHMIFteps}. Coming back to (\ref{EsupXi}), we finally observe that 
\begin{align*} 
\mathbb E\Big[\sup_{[0,T]}|X_t^{1,N,\epsilon}|\Big]\leq C\Big(\mathbb E[|X_0^1|]+T+\int_0^TI(F_s^{N,\epsilon})ds\Big)\leq C(1+\mathbb E[|X_0^1|]+T),
\end{align*}
which gives \eqref{EsupXi2} and concludes the proof.\hfill$\square$
\vskip0.5cm

We can now give the proof of existence and uniqueness for the particle system \eqref{systpart}.

\textbf{Proof of Theorem \ref{HMIFt}.} 
Like in \cite{TAK}, the key point of the proof is to show that particles of the system \eqref{systpart} a.s. never collide. We divide the proof in three steps. The first step consists in showing that a.s. there are no collisions between particles for the system \eqref{systparteps}. In the second step, we deduce that the particles of the system \eqref{systpart} also never collide, which ensures global existence and uniqueness for \eqref{systpart}. In the last step, we establish the estimates about the entropy, Fisher information and the first moment. We fix $N\geq2$ and for all $\epsilon\in(0,1)$, we consider $(X_t^{i,N,\epsilon})_{i=1,...,N,t\geq0}$ the unique solution to \eqref{systparteps}.

\textit{Step 1.} Let $\tau_\epsilon:=\inf\{t\geq0, \exists i\neq j, |X_t^{i,N,\epsilon}-X_t^{j,N,\epsilon}|\leq\epsilon\}$. The aim of this step is to prove that $\lim_{\epsilon\rightarrow0}\mathbb P[\tau_\epsilon<T]=0$ for all $T>0$. We fix $T>0$ and introduce
\begin{align} \label{Steps}
S_t^\epsilon:=\frac{1}{N^2}\sum_{i\neq j}\log|X_t^{i,N,\epsilon}-X_t^{j,N,\epsilon}|.
\end{align}
For any $A>1$, we have
\begin{align} \label{probataueps}
\mathbb P[\tau_\epsilon<T]&\leq\mathbb P\Big[\inf_{[0,T]}S_{t\wedge\tau_\epsilon}^\epsilon\leq S_{\tau_\epsilon}^\epsilon\Big]\\
\nonumber
&\leq \mathbb P[\exists i,\exists t\in[0,T],|X_{t}^{i,N,\epsilon}|>A]\\
\nonumber
&\ \ +\mathbb P\Big[\forall i, \forall t\in[0,T], |X_{t}^{i,N,\epsilon}|\leq A, \inf_{[0,T]}S_{t\wedge\tau_\epsilon}^\epsilon\leq S_{\tau_\epsilon}^\epsilon\Big]\\
\nonumber
&\leq \frac{N\mathbb E\Big[\sup_{[0,T]}|X_t^{1,N,\epsilon}|\Big]}{A}+\mathbb P\Big[\inf_{[0,T]}S_{t\wedge\tau_\epsilon}^\epsilon \leq \frac{\log\epsilon}{N^2}+\log 2A\Big]\\
\nonumber
&\leq \frac{C(1+T)N}{A}+\mathbb P\Big[\inf_{[0,T]}S_{t\wedge\tau_\epsilon}^\epsilon \leq \frac{\log\epsilon}{N^2}+\log 2A\Big],
\end{align}
where we used \eqref{EsupXi2}. We thus want to compute $\mathbb P\Big[\inf_{[0,T]}S_{t\wedge\tau_\epsilon}^\epsilon \leq-M\Big]$ for all (large) $M>0$. Using It\^o's formula, that $K_\epsilon(x)=K(x)$ for any $|x|\geq\epsilon$ (see \eqref{Keps}) and that $\triangle(\log|x|)=0$ on 
$\{x\in\mathbb R^2, |x|>\epsilon\}$, we have 
\begin{align*}
\log|X_{t\wedge\tau_\epsilon}^{i,N,\epsilon}&-X_{t\wedge\tau_\epsilon}^{j,N,\epsilon}|=\log |X_0^{i,N}-X_0^{j,N}|+M_{t\wedge\tau_\epsilon}^{i,j,\epsilon}\\
&\quad-\frac{\chi}{N}\int_0^{t\wedge\tau_\epsilon}\Big[\sum_{k\neq i,j}\Big(K(X_s^{i,N,\epsilon}-X_s^{k,N,\epsilon})-K(X_s^{j,N,\epsilon}-X_s^{k,N,\epsilon})\Big)\\
&\qquad\qquad\qquad+2K(X_s^{i,N,\epsilon}-X_s^{j,N,\epsilon})\Big].\frac{X_s^{i,N,\epsilon}-X_s^{j,N,\epsilon}}{|X_s^{i,N,\epsilon}-X_s^{j,N,\epsilon}|^2} ds\\
&\qquad\qquad\ \ =: \log |X_0^{i,N}-X_0^{j,N}|+M_{t\wedge\tau_\epsilon}^{i,j,\epsilon}+R_{t\wedge\tau_\epsilon}^{i,j,\epsilon},
\end{align*}
where $M_t^{i,j,\epsilon}$ is a martingale. Setting $S_0:=\frac{1}{N^2}\sum_{i\neq j}\log |X_0^{i,N}-X_0^{j,N}|$, $M_t^\epsilon:=\frac{1}{N^2}\sum_{i\neq j}M_{t\wedge\tau_\epsilon}^{i,j,\epsilon}$ and $R_t^\epsilon:=\frac{1}{N^2}\sum_{i\neq j}R_{t\wedge\tau_\epsilon}^{i,j,\epsilon}$, we thus have 
\begin{align*}
S_{t\wedge\tau_\epsilon}^\epsilon=S_0+M_t^\epsilon+R_t^\epsilon,
\end{align*}
so that 
\begin{align} \label{infSteps}
\mathbb P(\inf_{[0,T]}S_{t\wedge\tau_\epsilon}^\epsilon\leq-M)&\leq \mathbb P(S_0\leq-M/3)+\mathbb P(\inf_{[0,T]}M_t^\epsilon\leq-M/3)\\
\nonumber
&\quad+\mathbb P(\inf_{[0,T]}R_t^\epsilon\leq-M/3).
\end{align} 
Using exchangeability and that $|K(x)|=|x|^{-\alpha}$, we clearly have for some constant $C$ independent of $N$ and $\epsilon$,
\begin{align}
\nonumber\mathbb E[\sup_{[0,T]}|R_t^\epsilon|]&\leq C\chi\int_0^T\mathbb E\Big[\frac{1}{|X_s^{1,N,\epsilon}-X_s^{2,N,\epsilon}|^{\alpha+1}}\Big]ds\\
\nonumber&\leq C\chi\int_0^T(1+I(F_{s2}^{N,\epsilon}))ds\\
\label{ineq: Rteps}&\leq C(1+T),
\end{align}
where we used Lemma \ref{lemFish}, the fact that $I(F_{t2}^{N,\epsilon})\leq I(F_t^{N,\epsilon})$ by Lemma \ref{lemFishmarg}, and finally Proposition \ref{HMIFteps}. We thus get
\begin{align} \label{infRteps}
\mathbb P(\inf_{[0,T]}R_t^\epsilon\leq-M/3)\leq \mathbb P(\sup_{[0,T]}|R_t^\epsilon|\geq M/3)\leq \frac{C(1+T)}{M}.
\end{align}
We now want to compute $\mathbb P(\inf_{[0,T]}M_t^\epsilon\leq-M/3)$. Using that $\log|x|\leq |x|$, we have 
\begin{align*}
S_t^\epsilon\leq\frac{1}{N^2}\sum_{i\neq j}(|X_t^{i,N,\epsilon}|+|X_t^{j,N,\epsilon}|)\leq \frac{2}{N}\sum_i |X_t^{i,N,\epsilon}|.
\end{align*}
Consequently,
\begin{align*}
M_t^\epsilon&\leq S_{t\wedge\tau_\epsilon}^\epsilon+\sup_{s\in[0,T]}|R_s^\epsilon|-S_0\\
&\leq \frac{2}{N}\sum_{i}\sup_{s\in[0,T]}|X_s^{i,N,\epsilon}|+\sup_{s\in[0,T]}|R_s^\epsilon|-S_0=:K^\epsilon-S_0=:Z^\epsilon.
\end{align*}
We have 
\begin{align}\label{probaMteps}
\mathbb P(\inf_{[0,T]}M_t^\epsilon\leq-M/3)\leq \mathbb P(Z^\epsilon\geq\sqrt{M/3})+\mathbb P(\inf_{[0,T]}M_t^\epsilon\leq-M/3,Z^\epsilon<\sqrt{M/3}).
\end{align}
Since $(M_t^\epsilon)_{t\geq0}$ is a continuous local martingale, there exists a Brownian Motion $\beta$ such that $M_t^\epsilon=\beta_{<M^\epsilon>_t}$. For $x\in\mathbb R$, we set $\sigma_x:=\inf\{t\geq0, \beta_t=x\}$. Using that $\sup_{[0,T]}M_t^\epsilon\leq Z^\epsilon$ a.s.,
\begin{align}
\nonumber
\mathbb P(\inf_{[0,T]}M_t^\epsilon\leq-M/3&,Z^\epsilon<\sqrt{M/3})\leq \mathbb P(\inf_{[0,T]}M_t^\epsilon\leq-M/3,\sup_{[0,T]}M_t^\epsilon<\sqrt{M/3})\\
\nonumber
&\leq \mathbb P(\sigma_{-M/3}\leq\sigma_{\sqrt{M/3}})\\
\label{probaMteps2}
&=\frac{\sqrt{M/3}}{M/3+\sqrt{M/3}}\leq\sqrt\frac{3}{M},
\end{align}
by classical results on the Brownian Motion.
Using (\ref{EsupXi2}) and \eqref{ineq: Rteps}, we get that $\mathbb E[K^\epsilon]\leq C(1+T)$ where $C$ does not depend on $\epsilon$. So using the Markov inequality, 
\begin{align}
\nonumber
\mathbb P(Z^\epsilon\geq \sqrt{M/3})&=\mathbb P(K^\epsilon-S_0\geq \sqrt{M/3})\\
\nonumber
&\leq \mathbb P(K^\epsilon\geq \sqrt{M/12})+\mathbb P(-S_0\geq \sqrt{M/12})\\
\label
{probaZeps}
&\leq \frac{C(1+T)}{\sqrt{M}}+\mathbb P(-S_0\geq \sqrt{M/12}).
\end{align}
Gathering \eqref{probaMteps}, \eqref{probaMteps2} and \eqref{probaZeps}, we find that
\begin{align} 
\label{infMteps}
\mathbb P(\inf_{[0,T]}M_t^\epsilon\leq-M/3)&\leq \frac{C(1+T)}{\sqrt{M}}+\mathbb P(-S_0\geq \sqrt{M/12}).
\end{align}
Coming back to \eqref{probataueps} and \eqref{infSteps}, using \eqref{infRteps} and \eqref{infMteps} with $M=-\frac{\log\epsilon}{N^2}-\log 2A$, we finally get
 that for any $\epsilon\in(0,1)$, any $A>1$ such that $\frac{\log\epsilon}{N^2}+\log 2A<0$,
\begin{align*}
\mathbb P(\tau_\epsilon<T)&\leq \frac{C(1+T)N}{A}+\mathbb P\Big(S_0\leq \big(\frac{\log\epsilon}{N^2}+\log 2A\big)/3\Big)\\
&\quad+\frac{C(1+T)}{-\frac{\log\epsilon}{N^2}-\log 2A}+\frac{C(1+T)}{\sqrt{\frac{-\log\epsilon}{N^2}-\log 2A}}\\
&\quad+\mathbb P\Big(S_0\leq -\sqrt{\big(-\frac{\log\epsilon}{N^2}-\log 2A\big)/12}\Big).
\end{align*}
Observe finally that $S_0>-\infty$ a.s. (because $F_0^N$ has a density since $H(F_0^N)<\infty$) so that $\lim_{M\rightarrow+\infty}\mathbb P(S_0<-M)=0$. Letting $\epsilon\rightarrow0$ in the above formula, we get that for all $A>1$,
\begin{align*}
\limsup_\epsilon\mathbb P(\tau_\epsilon<T)\leq \frac{C(1+T)N}{A}.
\end{align*}
It only remains to make $A$ go to $\infty$ to conclude this step. 

\textit{Step 2.} Since $K$ is Lipschitz-continuous outside $0$, classical arguments give existence and uniqueness of a solution to \eqref{systpart} until the explosion time $\tau=\inf\{t\geq0, \exists i\neq j, X_t^{i,N}=X_t^{j,N}\}$. We can observe that since $K_\epsilon(x)=K(x)$ for any $|x|\geq\epsilon$, $(X^{i,N,\epsilon})_{
i=1,...,N}$ is solution to \eqref{systpart} on $[0,\tau_\epsilon]$ so that for any $i=1,...,N$, $X_t^{i,N}=X_t^{i,N,\epsilon}$ on $[0,\tau_\epsilon]$. We thus have $\tau_\epsilon<\tau$ for any $\epsilon\in(0,1)$ a.s. so that, using Step 1, we have for any $T>0$
\begin{align*}
\mathbb P(\tau<T)\leq \mathbb P(\tau_\epsilon<T)\underset{\epsilon\rightarrow0}{\longrightarrow }0.
\end{align*} 
Thus $\tau=\infty$ a.s. which proves global existence and uniqueness for \eqref{systpart}.

\textit{Step 3.} Using that the functionals $H$, $I$ and $M_1$ are lower semi-continuous and Proposition \ref{HMIFteps}, we have
\begin{align}
\nonumber
H(F_t^N)&\leq\liminf_\epsilon H(F_t^{N,\epsilon})\leq C(1+t),\\
\label{I(Fsn)demo}
\int_0^tI(F_s^N)ds&\leq\liminf_\epsilon\int_0^t I(F_s^{N,\epsilon})ds\leq C(1+t),
\end{align}
and
\begin{align*}
M_1(F_t^{N})\leq  \liminf_\epsilon M_1(F_t^{N,\epsilon})\leq C(1+t).
\end{align*}
Using the Fatou lemma and \eqref{EsupXi2}, we get
\begin{align*}
\mathbb E\Big[\sup_{[0,T]}|X_t^{1,N}|\Big]\leq \liminf_\epsilon\mathbb E\Big[\sup_{[0,T]}|X_t^{1,N,\epsilon}|\Big]\leq C(1+T),
\end{align*}
and \eqref{EsupX1} is proved. It remains to prove \eqref{entroFishFt}. Using again that the functionals $H$ and $I$ are lower semi-continuous and using \eqref{entroFteps}, we get
\begin{align*}
H(F_t^N)+\int_0^tI(F_s^N)ds&\leq \liminf_{\epsilon}\Big[H(F_t^{N,\epsilon})+\int_0^tI(F_s^{N,\epsilon})ds\Big]\\
&\leq H(F_0^N)+\liminf_\epsilon \frac{\chi}{N^2}\int_0^t \sum_{i\neq j}\mathbb E[{\rm div}K_\epsilon(X_s^{i,N,
\epsilon}-X_s^{j,N,\epsilon})]ds.
\end{align*}
By exchangeability, it suffices to prove that, as $\epsilon\rightarrow0$,
\begin{align*}
D_\epsilon:=\int_0^t\mathbb E[{\rm div}K_\epsilon(X_s^{1,N,\epsilon}-X_s^{2,N,\epsilon})]ds\rightarrow\int_0^t\mathbb E[{\rm div}K(X_s^{1,N}-X_s^{2,N})]ds=:D.
\end{align*}
By Step 2, we have $X_s^{i,N}=X_s^{i,N,\epsilon}$ for any $i$ and $s\leq\tau_\epsilon$ and thus recalling that $K_\epsilon(x)=K(x)$ for any $|x|\geq\epsilon$, we get that a.s. for any $s<\tau_\epsilon$
\begin{align*}
{\rm div}K_\epsilon(X_s^{1,N,\epsilon}-X_s^{2,N,\epsilon})={\rm div}K(X_s^{1,N,\epsilon}-X_s^{2,N,\epsilon})={\rm div}K(X_s^{1,N}-X_s^{2,N}).
\end{align*}
So using that div$K(x)\leq2|x|^{-\alpha-1}$ and div$K_\epsilon(x)\leq2|x|^{-\alpha-1}$, we get
\begin{align*}
|D-D_\epsilon|\leq C\int_0^t \mathbb E\Big[\one_{\{\tau_\epsilon<s\}}\Big(\frac{1}{|X_s^{1,N,\epsilon}-X_s^{2,N,\epsilon}|^{\alpha+1}}+\frac{1}{|X_s^{1,N}-X_s^{2,N}|^{\alpha+1}}\Big)\Big]ds.
\end{align*}
Let $a\in\Big(0,\frac{1-\alpha}{1+\alpha}\Big)$ (in order to have $(1+a)(\alpha+1)<2$). Using first the H\"older inequality with $p=1+a$ and $q$ such that $1/p+1/q=1$, and then Lemma \ref{lemFish} with $\beta=1$, we get 
\begin{align*}
|D-D_\epsilon|&\leq C\int_0^t \mathbb P(\tau_\epsilon<s)^{1/q}\mathbb E\Big[\Big(\frac{1}{|X_s^{1,N,\epsilon}-X_s^{2,N,\epsilon}|^{(\alpha+1)(1+a)}}\\
&\qquad\qquad\qquad\qquad\qquad+\frac{1}{|X_s^{1,N}-X_s^{2,N}|^{(\alpha+1)(1+a)}}\Big)\Big]^{1/p}ds\\
&\leq C \mathbb P(\tau_\epsilon<t)^{1/q}\int_0^t[1+I(F_s^{N,\epsilon})+I(F_s^N)]ds\\
&\leq C(1+t) \mathbb P(\tau_\epsilon<t)^{1/q},
\end{align*}
by \eqref{eqHMIFteps} and \eqref{I(Fsn)demo}. This tends to $0$ as $\epsilon\rightarrow0$ by Step1 and concludes the proof. \hfill$\square$
 
\section{Convergence of the particle system}

We start this section with a tightness result for the particle system \eqref{systpart}.
\begin{lemme} \label{tension}
Let $\alpha\in(0,1)$. Assume \eqref{f0chaotic}. For each $N\geq2$, 
let $(X_t^{i,N})_{i=1,...,N}$ be the unique solution to \eqref{systpart} and $Q^N:=\frac{1}{N}\sum_{i=1}^N\delta_{(X_t^{i,N})_{t\geq0}}$.\\
(i) The family $\{\mathcal{L}((X_t^{1,N})_{t\geq0}),N\geq2\}$ is tight in $\textbf{P}(C([0,\infty),\mathbb{R}^2))$.\\
(ii) The family $\{\mathcal{L}(Q^N),N \geq2\}$ is tight in $\textbf{P}(\textbf{P}(C([0,\infty),\mathbb{R}^2)))$.
\end{lemme}

\textbf{Proof.}
Since the system is exchangeable, we deduce $(ii)$ from $(i)$ by \cite[Proposition 2.2]{SZN}. Let's prove $(i)$. Let thus $\eta>0$ and $T>0$ be fixed. To prove the tightness of $\{\mathcal{L}((X_t^{1,N})_{t\geq0}),N\geq2\}$ in $\textbf{P}(C([0,\infty),\mathbb{R}^2))$, we have to find a compact subset $\mathcal{K}_{\eta,T}$ of $C([0,T],\mathbb{R}^2)$ such that $\sup_N \mathbb{P}[(X_t^{1,N})_{t\in[0,T]})\notin\mathcal{K}_{\eta,T}]\leq\eta$. We first set $Z_T:=\sup_{0<s<t<T}\sqrt{2}|B_t^1-B_s^1|/|t-s|^{1/3}$. This random variable is a.s. finite since the paths of a Brownian motion are a.s. H\"older continuous with index $1/3$. We can also notice that the law of $Z_T$ does not depend on $N$. Using the H\"older inequality with $p=3$ and $q=3/2$, we get that for all $0<s<t<T$,
\begin{align*}
\Big|\frac{\chi}{N}\sum_{j=2}^N\int_s^tK(X_u^{1,N}-X_u^{j,N})du\Big|&\leq \frac{\chi}{N}\sum_{j=2}^N\int_s^t\frac{du}{|X_u^{1,N}-X_u^{j,N}|^{\alpha}}\\
&\leq \frac{\chi}{N}(t-s)^{1/3}\sum_{j=2}^N\Big(\int_0^T\frac{du}{|X_u^{1,N}-X_u^{j,N}|^{3\alpha/2}}\Big)^{2/3}\\
&\leq (t-s)^{1/3} \Big(\chi+\frac{\chi}{N}\sum_{j=2}^N\int_0^T\frac{du}{|X_u^{1,N}-X_u^{j,N}|^{3\alpha/2}}\Big)\\
&=:(t-s)^{1/3}U_T^N.
\end{align*}   
Using Lemma \ref{lemFish} with $\gamma=3\alpha/2$ and $\beta=1$, the exchangeability of the system of particles, and denoting by $F_{u2}^N$ the two-marginal of $F_u^N$, we have
\begin{align*}
\mathbb{E}(U_T^N)=\chi+\chi\frac{N-1}{N}\int_0^T \mathbb{E}\Big(\frac{1}{|X_u^{1,N}-X_u^{2,N}|^{3\alpha/2}}\Big)du&\leq \chi+C\int_0^T(1+I(F_{u2}^N))du\\
&\leq \chi+C\int_0^T(1+I(F_u^N))du\\
&\leq C(1+T),
\end{align*}
where we used that $I(F_{t2}^N)\leq I(F_t^{N,\epsilon})$ by Lemma \ref{lemFishmarg} and Theorem \ref{HMIFt}. We thus have $\sup_{N\geq2}\mathbb{E}(U_T^N)<\infty$. Furthermore, $Z_T$ is also a.s. finite so that we can find $R>0$ such that $\mathbb{P}(Z_T+U_T^N>R)\leq\eta/2$ for all $N\geq2$. Recalling (\ref{f0chaotic}), we can also find $a>0$ such that $\sup_{N\geq2}\mathbb{P}(X_0^{1,N}>a)\leq\eta/2$. We now consider
\begin{align*}
\mathcal{K}_{\eta,T}:=\{f\in C([0,T],\mathbb R^2), |f(0)|\leq a, |f(t)-f(s)|\leq R(t-s)^{1/3}\ \forall 0<s<t<T\},
\end{align*}
which is a compact subset of $C([0,T],\mathbb R^2)$ by Ascoli's theorem. Observing that for all $0<s<t<T$,
$|X_t^{1,N}-X_s^{1,N}|\leq(Z_T+U_t^N)(t-s)^{1/3}$, we get
\begin{align*}
\mathbb{P}[(X_t^{1,N})_{t\in[0,T]}\notin\mathcal{K}_{\eta,T}]\leq \mathbb{P}(|X_0^{1,N}|>a)+\mathbb{P}(Z_T+U_T^N>R)\leq\eta,
\end{align*}
which concludes the proof. \hfill$\square$
\vskip0.5cm

We define $\mathcal S$ as the set of all probability measures $f\in\textbf{P}(C([0,\infty),\mathbb{R}^2))$ such that $f$ is the law of $(X_t)_{t\geq0}$ solution to (\ref{edsKS}) satisfying (setting $f_t=\mathcal L(X_t)$)
\begin{align}\label{eq:FisherInf}
\forall T>0,\quad \int_0^T I(f_s)ds<\infty\quad{\rm and}\quad \sup_{[0,T]}M_1(f_s)<\infty.
\end{align}
Observe that by Lemma \ref{lemLpFish}, \eqref{eq:FisherInf} implies \eqref{eq:FinFirstMom}.

\begin{prop} 
\label{CVversX}
Let $\alpha\in(0,1)$ and assume \eqref{f0chaotic}. For each $N\geq2$, let $(X_0^{i,N})_{i=1,...,N}$ be $F_0^N$-distributed and consider the solution $(X_t^{i,N})_{i=1,...,N, t\geq0}$ to (\ref{systpart}). Assume that there is a subsequence of $\mathcal{Q}^N:=\frac{1}{N}\sum_{i=1}^N\delta_{(X_t^{i,N})_{t\geq0}}$ going in law to some $\textbf{P}(C([0,\infty),\mathbb{R}^2))$-valued random variable $\mathcal{Q}$. Then $\mathcal{Q}$ a.s. belongs to $\mathcal S$.
\end{prop}

\textbf{Proof.} We consider a (not relabelled) subsequence of $\mathcal Q^N$ going in law to some $\mathcal Q$ and we introduce the identity map $\psi:C([0,\infty);\mathbb R^2)\rightarrow C([0,\infty);\mathbb R^2)$. Using the arguments of \cite[Proposition 6.1]{FHM}, we have to prove that $\mathcal Q$ a.s. satisfies
\begin{itemize}
\item[(a)] $\mathcal Q\circ(\psi(0))^{-1}=f_0$;
\item[(b)] setting $\mathcal Q_t=\mathcal Q\circ(\psi(t))^{-1}$, $(\mathcal Q_t)_{t\geq0}$ satisfies \eqref{eq:FisherInf};
\item[(c)] for all $0<t_0<\ldots<t_k<s<t$, $\varphi_1,\ldots,\varphi_k\in C_b(\mathbb R^2)$, $\varphi\in C_b^2(\mathbb R^2)$, $\mathcal F(\mathcal Q)=0$ where, for $f\in\textbf{P}(C([0,\infty),\mathbb{R}^2))$,
\begin{align*}
\mathcal F(f):=&\iint f(d\gamma)f(d\tilde\gamma)\varphi_1(\gamma_{t_1})\ldots\varphi_k(\gamma_{t_k})\\ 
&\left[\varphi(\gamma_t)-\varphi(\gamma_s)+\chi\int_s^t\nabla_x\varphi(\gamma_u)\cdot K(\gamma_u-\tilde\gamma_u)\,du-\int_s^t\triangle_x\varphi(\gamma_u)du\right].
\end{align*}
\end{itemize}
For simplicity, we split the proof in many steps.

\textit{Step 1.} By assumption \eqref{f0chaotic}, we have that $F^N_0$ is $f_0$-chaotic which implies that $\mathcal{Q}_0^N=\mathcal{Q}^N\circ\psi(0)^{-1}$ goes weakly to $f_0$ in law, and, since $f_0$ is deterministic, also in probability. Hence $\mathcal{Q}_0=f_0$ a.s. and thus $f\circ\psi(0)^{-1}=f_0$. Thus $\mathcal Q$ a.s. satisfies (a).\\

\textit{Step 2.} Since $\frac{1}{N}\sum_{i=1}^N\delta_{X_t^{i,N}}$ goes weakly to $\mathcal Q_t$, for all $j\geq 1, F_{tj}^N$ goes weakly to $\pi_{tj}$, where $\pi_t:=\mathcal L(\mathcal Q_t)$ and $\pi_{tj}:=\int_{\textbf{P}(\mathbb R^2)}f
^{\otimes j}\pi_t(df)$. We can thus apply \cite[Theorem 5.7]{HM} (and then the Fatou Lemma) to get
\begin{eqnarray*}
 \mathbb E\Big[\int_0^TI(\mathcal Q_s)ds\Big]=\int_0^T \mathbb E[I(\mathcal Q_s)]ds&\leq&\int_0^T\liminf_N I(F^N_s)\,ds\\&\leq&\liminf_N \int_0^TI(F^N_s)\,ds,
\end{eqnarray*}
which is finite by Theorem \ref{HMIFt}. We conclude that $\int_0^TI(\mathcal Q_s)ds<\infty$ a.s. We also have, using the Fatou lemma and the exchangeability of the particles,
\begin{align*}
\mathbb E\Big[\sup_{[0,T]}M_1(\mathcal Q_t)\Big]&\leq\mathbb E\Big[\liminf_N\sup_{[0,T]}M_1(\mathcal Q_t^N)\Big]\\
&\leq \liminf_N\mathbb E\Big[\sup_{[0,T]}\frac{1}{N}\sum_{i=1}^N|X_t^{i,N}|\Big]\\
&\leq \liminf_N\mathbb E\Big[\sup_{[0,T]}|X_t^{1,N}|\Big]\leq C(1+T),
\end{align*}
by \eqref{EsupX1}, so that $\sup_{[0,T]}M_1(\mathcal Q_t)<\infty$ a.s. Consequently, $\mathcal Q$ a.s. satisfies (b).\\

\textit{Step 3.1.} 
Using It\^o's formula
\begin{align*}
 O_t^i:=&\varphi(X_s^{i,N})+\frac{\chi}{N}\sum_{j\neq i}\int_0^t\nabla_x\varphi(X_s^{i,N}))\cdot K(X_s^{i,N}-X_s^{j,N}))ds-\int_0^t\triangle_x\varphi(X_s^{i,N})ds\\
 =&\varphi(X_0^{i,N})+\sqrt{2}\int_0^t\nabla_x\varphi(X_s^{i,N})\cdot dB_s^i.
\end{align*}
But, using the last equality, we see that
\begin{align*}
\mathcal F(\mathcal Q^N)&=\frac1N\sum_{i=1}^N\varphi_1(X_{t_1}^{i,N})\ldots\varphi_k(X_{t_k}^{i,N})[O_t^i-O_s^i]\\
&=\frac{\sqrt{2}}N\sum_{i=1}^N\varphi_1(X_{t_1}^{i,N})\ldots
\varphi_k(X_{t_k}^{i,N})\int_s^t\nabla_x\varphi(X_u^{i,N})\cdot dB_u^i.
\end{align*}
From there, and thanks to the independence of the Brownian motions we conclude that (recall that the functions $\varphi_1,...,\varphi_k,\nabla_x\varphi$ are bounded) $$\mathbb E\left[(\mathcal F(\mathcal Q^N))^2\right]\leq \frac{C_{\mathcal F}}{N}.$$

\textit{Step 3.2.} We also introduce the regularized version of $\mathcal F$. For $\varepsilon\in(0,1)$, we define $\mathcal F_\varepsilon$ replacing $K$ by $K_\varepsilon$ defined by \eqref{Keps}. Since $f\mapsto\mathcal F_\varepsilon(f)$ is continuous and bounded from ${\bf P}(C([0,\infty);\mathbb R^2))$ to $\mathbb R$ and since $\mathcal Q^N$ goes in law to $\mathcal Q$, we deduce that for any $\varepsilon\in(0,1)$, $$\mathbb E\mathbb[|\mathcal F_\varepsilon(
\mathcal Q)|]=\lim_N\mathbb E[|\mathcal F_\varepsilon(\mathcal Q^N)|].$$

\textit{Step 3.3.} Using that all the functions and their derivatives involved in $\mathcal F$ are bounded and that $|K_\varepsilon(x)-K(x)|\leq |x|^{-\alpha}\one_{0\leq|x|\leq\varepsilon}$, we get
\begin{align*}
 |\mathcal F(f)-\mathcal F_\varepsilon(f)|\leq&\,\chi\,C_{\mathcal F}\iiint_0^t|\gamma(u)-\tilde\gamma(u)|^{-
\alpha}\one_{0<|\gamma(u)-\tilde\gamma(u)|<\varepsilon}\,du\,f(d\gamma)f(d\tilde\gamma)\\
 \leq&\, C_{\mathcal F}\varepsilon^{3/2-\alpha}\iiint_0^t|\gamma(u)-\tilde\gamma(u)|^{-3/2}\one_{\gamma(u)\neq\tilde\gamma(u)}du\,f(d\gamma)f(d\tilde\gamma).
\end{align*}
Thus,
\begin{equation}\nonumber
 |\mathcal F(\mathcal Q^N)-\mathcal F_\varepsilon(\mathcal Q^N)|\leq\frac{C_{\mathcal F}\varepsilon^{3/2-\alpha}}{N^2}\sum_{i\neq j}\int_0^t|X_u^{i,N}-X_u^{j,N}|^{-3/2}\,du,
\end{equation}
and by exchangeability
\begin{equation}\nonumber
 \mathbb E\left[|\mathcal F(\mathcal Q^N)-\mathcal F_\varepsilon(\mathcal Q^N)|\right]\leq C_{\mathcal F}\varepsilon^{3/2-\alpha}\int_0^t\mathbb E\left[|X_u^{1,N}-X_u^{2,N}|^{-3/2}\right]\,du.
\end{equation}
Using Lemma \ref{lemFish} with $\gamma=3/2$ and $\beta=1$ and denoting by $F_{u2}^N$ the two-marginal of $F_u^N$, we have
\begin{equation}\nonumber
 \mathbb E\left[|\mathcal F(\mathcal Q^N)-\mathcal F_\varepsilon(\mathcal Q^N)|\right]\leq C_{\mathcal F}\varepsilon^{3/2-\alpha}\int_0^tI(F^N_{u2})\,du.
\end{equation}
Using that $I(F_{t2}^N)\leq I(F_t^{N})$ by Lemma \ref{lemFishmarg} and Theorem \ref{HMIFt} we conclude that
 \begin{equation}\nonumber
 \mathbb E\left[|\mathcal F(\mathcal Q^N)-\mathcal F_\varepsilon(\mathcal Q^N)|\right]\leq C_{\mathcal F}\varepsilon^{3/2-\alpha}.
\end{equation}

\textit{Step 3.4.} Now we see that
\begin{align*}
 |\mathcal F(\mathcal Q)-\mathcal F_\varepsilon(\mathcal Q)|\leq&\, C_{\mathcal F}\varepsilon^{3/2-\alpha}\int_0^t\int_
{\mathbb R^2}\int_{\mathbb R^2}|x-y|^{-3/2}\mathcal Q_s(dx)\mathcal Q_s(dy)\,ds.
\end{align*}
Step 2 says that \eqref{eq:FisherInf} holds true for $\mathcal Q_s$, then thanks to Lemma \ref{lemLpFish} we get that a.s., $\nabla_x\mathcal Q_s\in L^{2q/(3q-2)}(0,T;L^q(\mathbb R^2))$ for all $q\in [1,2)$. Then using \cite[Lemma 3.5]{FHM} for $\gamma=3/2$ we deduce that a.s. $$ \lim_{\varepsilon\rightarrow0}|\mathcal F(\mathcal Q)-\mathcal F_\varepsilon(\mathcal Q)|= 0.$$

\textit{Step 3.5.} Using Steps 3.1, 3.2 and 3.3, we finally observe, using the same arguments as in \cite[Proposition 6.1, \textit{Step 4.5}]{FHM}, that
\begin{align*}
\mathbb E[|\mathcal F(\mathcal Q)|\wedge1] \leq C_{\mathcal F}\varepsilon^{3/2-
\alpha}+\mathbb E[|\mathcal F(\mathcal Q)-\mathcal F_\varepsilon(\mathcal Q)|\wedge1],
\end{align*}
so that $\mathcal F(\mathcal Q)=0$ a.s. by Step 3.4 thanks to dominated convergence and $\mathcal Q$ a.s. satisfies (c) which concludes the proof. \hfill$\square$

\section{Well-posedness and propagation of chaos}

We start this section with the proof of existence and uniqueness for the nonlinear S.D.E. \eqref{edsKS}. We will use that for $\gamma\in(-2,0)$, for $q\in(2/(2+\gamma),\infty]$ and for any $h\in\textbf{P}(\mathbb R^2)\cap L^q(\mathbb R^2)$,
\begin{align}
\nonumber
\sup_{v\in\mathbb R^2}\int_{\mathbb R^2}h(v_*)|v-v_*|^\gamma dv_*&\leq \sup_{v\in\mathbb R^2}\int_{|v_*-v|<1}h(v_*)|v-v_*|^\gamma dv_*\\
\nonumber
&\quad+\sup_{v\in\mathbb R^2}\int_{|v_*-v|\geq1}h(v_*)dv_*\\
\label{Jalpha}
&\leq C_{\gamma,q}||h||_{L^q(\mathbb R^2)}+1,
\end{align}
where 
\begin{align*}
C_{\gamma,q}=\Big[\int_{|v_*|\leq1}|v_*|^{\gamma q/(q-1)}dv_*\Big]^{(q-1)/q}<\infty,
\end{align*}
since by assumption $\gamma q/(q-1)>-2$.

\textbf{Proof of Theorem \ref{exuniedsKS}.}
The existence in law follows from Proposition \ref{CVversX} and Lemma \ref{tension} (see the comment after \eqref{eq:FisherInf}). We now prove pathwise uniqueness which will also imply the strong existence. To this aim, we consider $(X_t)_{t\geq0}$ and $(Y_t)_{t\geq0}$ two solutions of (\ref{edsKS}) such that, setting $f_s:=\mathcal L(X_s)$ and $g_s:=\mathcal L(Y_s)$,  $(f_t)_{t\geq0}$ and $(g_t)_{t\geq0}$ are in $L_{loc}^\infty([0,\infty),\textbf{P}_1(\mathbb R^2))\cap L_{loc}^1([0,\infty);L^p(\mathbb 
R^2))$ for some $p>\frac{2}{1-\alpha}$. For any $s>0$, we consider the probability measure $R_s$ on $\mathbb R^2\times\mathbb R^2$ with first (respectively second) marginal equal to $f_s$ (resp. $g_s$) such that 
\begin{align*}
\mathcal W_1(f_s,g_s)=\int_{\mathbb R^2\times\mathbb R^2}|x-y|R_s(dx,dy).
\end{align*}
We have
\begin{align*}
X_t-Y_t&=-\chi\Big(\int_0^t\int_{\mathbb R^2}K(X_s-x)f_s(dx)ds-\int_0^t\int_{\mathbb R^2}K(Y_s-y)g_s(dy)ds\Big)\\
&=-\chi\int_0^t\int_{\mathbb R^2\times\mathbb R^2}[K(X_s-x)-K(Y_s-x)]R_s(dx,dy).
\end{align*}
Using Lemma \ref{Kx-Ky} and recalling that $\mathcal L(X_t)=f_t$, $\mathcal L(Y_t)=g_t$, and that $R_t$ has marginals $f_t$ and $g_t$, this gives 
\begin{align*}
\mathbb E[\sup_{[0,T]}|X_t-Y_t|]&\leq C_\alpha\chi\int_0^T\int_{\mathbb R^2\times\mathbb R^2}\mathbb E\Big[(|X_s-Y_s|+|x-y|)\Big(\frac{1}{|X_s-x|^{\alpha+1}}\\
&\qquad\qquad\qquad\qquad\qquad\qquad\qquad\quad+\frac{1}{|Y_s-y|^{\alpha+1}}\Big)\Big]R_s(dx,dy)ds\\
&\leq C_\alpha\chi\int_0^T\mathbb E\Big[|X_s-Y_s|\Big(\int_{\mathbb R^2}\frac{1}{|X_s-x|^{\alpha+1}}f_s(dx)\\
&\qquad\qquad\qquad\qquad\qquad\qquad+\int_{\mathbb R^2}\frac{1}{|Y_s-y|^{\alpha+1}}g_s(dy)\Big)\Big]ds\\
&\quad +C_\alpha\chi\int_0^T\int_{\mathbb R^2\times\mathbb R^2}|x-y|\mathbb E\Big[\frac{1}{|X_s-x|^{\alpha+1}}\\
&\qquad\qquad\qquad\qquad\qquad\qquad\qquad+\frac{1}{|Y_s-y|^{\alpha+1}}\Big]R_s(dx,dy)ds.
\end{align*}
Using \eqref{Jalpha}, we thus have, since $\int_{\mathbb R^2\times\mathbb R^2}|x-y|R_s(dx,dy)=\mathcal W_1(f_s,g_s)\leq \mathbb E[|X_s-Y_s|]$ by definition of $\mathcal W_1$,
\begin{align*}
\mathbb E[\sup_{[0,T]}|X_t-Y_t|]&\leq C\int_0^T\mathbb E[|X_s-Y_s|](1+||f_s||_{L^p}+||g_s||_{L^p})ds\\
&\quad+C\int_0^T\int_{\mathbb R^2\times\mathbb R^2}|x-y|(1+||f_s||_{L^p}+||g_s||_{L^p})R_s(dx,dy)ds\\
&\leq C\int_0^T\mathbb E[|X_s-Y_s|](1+||f_s+g_s||_{L^p})ds.
\end{align*}
By Gr\"onwall's Lemma, we thus get $\mathbb E(\sup_{[0,T]}|X_t-Y_t|)=0$ and pathwise uniqueness is proved. \hfill$\square$
\vskip 0.5cm

The following lemma is useful for the uniqueness of \eqref{eq:KS}.

\begin{lemme} \label{lemuni}
Let $p>2/(1-\alpha)$ and consider a weak solution $(f_t)_{t\geq0}$ to \eqref{eq:KS} lying in $L_{loc}^\infty([0,\infty),\textbf{P}_1(\mathbb R^2))\cap L_{loc}^1([0,\infty);L^p(\mathbb R^2))$. Assume that for some $h=(h_t)_{t\geq0}$ lying in $L_{loc}^\infty([0,\infty),\textbf{P}_1(\mathbb R^2))\cap L_{loc}^1([0,\infty);L^p(\mathbb R^2))$, for all $\varphi\in C_c^2(\mathbb R^2)$, all $t\geq0$,
\begin{align} \label{KSlinweak}
\int_{\mathbb R^2}\varphi(x)h_t(dx)=&\int_{\mathbb R^2}\varphi(x)f_0(dx)+\int_0^t\int_{\mathbb R^2}\triangle_x\varphi(x)h_s(dx)\,ds\\
\nonumber
 &-\chi\int_0^t\int_{\mathbb R^2}\int_{\mathbb R^2}K(x-y)\cdot\nabla_x\varphi(x)f_s(dy)h_s(dx)\,ds.
\end{align}
Then $h=f$.
\end{lemme}

\textbf{Proof.} For any $\varphi\in C_c^2(\mathbb R^2)$ and any $t\geq0$, we set
\begin{align*}
\mathcal A_t\varphi(x)=\triangle_x\varphi(x)-\chi\int_{\mathbb R^2}K(x-y)\cdot\nabla_x\varphi(x)f_t(dy).
\end{align*}
We will prove that for any $\mu\in\textbf{P}_1(\mathbb R^2)$, there exists at most one $h$ lying in $L_{loc}^\infty([0,\infty),\textbf{P}_1(\mathbb R^2))\cap L_{loc}^1([0,\infty);L^p(\mathbb R^2))$ such that for all $t\geq0$, $\varphi\in C_c^2(\mathbb R^2)$,
\begin{align} \label{KSlin}
\int_{\mathbb R^2}\varphi(x)h_t(dx)=\int_{\mathbb R^2}\varphi(x)\mu(dx)+\int_0^t\int_{\mathbb R^2}\mathcal A_s\varphi(x)h_s(dx)ds.
\end{align}
This will conclude the proof since $f$ and $h$ solve this equation with $\mu=f_0$ by assumption.

\textit{Step 1.} Let $\mu\in\textbf{P}_1(\mathbb R^2)$. A continuous adapted $\mathbb R^2$-valued process $(X_t)_{t\geq0}$ on some filtered probability space $(\Omega,\mathcal F,(\mathcal F_t)_{t\geq0},P)$ is said to solve the martingale problem $MP((\mathcal A_t)_{\geq0},\mu)$ if $P\circ X_0^{-1}=\mu$ and if for all $\varphi\in C_c^2(\mathbb R^2)$, $(M_t^\varphi)_{t\geq0}$ is a $(\Omega,\mathcal F,(\mathcal F_t)_{t\geq0},P)$-martingale, where
\begin{align*}
M_t^\varphi=\varphi(X_t)-\int_0^t\mathcal A_s\varphi(X_s)ds.
\end{align*}
Using Bhatt-Karandikar \cite[Theorem 5.2]{BK} (see also Remark 3.1 in \cite{BK}), uniqueness for \eqref{KSlin} holds if\\
(i) there exists a countable subset $(\varphi_k)_{k\geq1}\subset C_c^2$ such that for all $t\geq0$, the closure (for the bounded pointwise convergence) of $\{(\varphi_k,\mathcal A_t\varphi_k),k\geq1\}$ contains $\{(\varphi,\mathcal A_t\varphi),\varphi\in C_c^2\}$,\\
(ii) for each $x_0\in\mathbb R^2$, there exists a solution to $MP((\mathcal A_t)_{\geq0},\delta_{x_0})$,\\
(iii) for each $x_0\in\mathbb R^2$, uniqueness (in law) holds for $MP((\mathcal A_t)_{\geq0},\delta_{x_0})$.
\vskip0.3cm

\textit{Step 2.} We first prove (i). Consider thus some countable $(\varphi_k)_{k\geq1}\subset C_c^2$ dense in $C_c^2$, in the sense that for $\psi\in C_c^2$, there exists a subsequence $\varphi_{k_n}$ such that $\lim_{n\rightarrow\infty}(||\psi-\varphi_{k_n}||_\infty+||\psi'-\varphi_{k_n}'||_\infty+||\psi''-\varphi_{k_n}''||_\infty)=0$. We then have to prove that, for $t
\geq0$,\\
(a) $\mathcal A_t\varphi_{k_n}(x)$ tends to $\mathcal A_t\psi(x)$ for all $x\in\mathbb R^2$,\\
(b)$\sup_n||\mathcal A_t\varphi_{k_n}||_{\infty}<\infty$.

Let $x\in\mathbb R^2$. By Lemma \ref{Kx-Ky}, we have
\begin{align*}
|\mathcal A_t\varphi_{k_n}(x)-\mathcal A_t\psi(x)|&\leq ||\psi''-\varphi_{k_n}''||_\infty+\chi||\psi'-\varphi_{k_n}'||_\infty\int_{\mathbb R^2}\frac{1}{|x-y|^\alpha}f_t(dy)\rightarrow0,
\end{align*}
since $\int_{\mathbb R^2}\frac{1}{|x-y|^\alpha}f_t(dy)\leq C(1+||f_t||_{L^p})$ by \eqref{Jalpha}. For (b), we can observe that setting $A:=\sup_n(||\varphi_{k_n}||_\infty+||\varphi_{k_n'}||_\infty+||\varphi_{k_n}''||_\infty)$
\begin{align*}
|\mathcal A_t\varphi_{k_n}|\leq A+\chi A\int_{\mathbb R^2}\frac{1}{|x-y|^\alpha}f_t(dy)\leq A+CA(1+||f_t||_{L^p}),
\end{align*}
which concludes this step.
\vskip0.3cm

\textit{Step 3.} 
Using classical arguments, we observe that a process $(X_t)_{t\geq0}$ is a solution to $MP((\mathcal A_t)_{\geq0},\delta_{x_0})$ if and only if there exists, on a possibly enlarged probability space, a $(\mathcal F_t)_{t\geq0}$-Brownian motion $(B_t)_{t\geq0}$ such that
\begin{align} \label{EDSKSlin}
X_t=x_0-\chi\int_0^t\int_{\mathbb R^2}K(X_s-x)f_s(dx)ds+\sqrt2 B_t.
\end{align}
It thus suffices to prove existence and uniqueness in law for solutions to \eqref{EDSKSlin} to get (ii) and (iii).

\textit{Step 4.} The proof of (pathwise) uniqueness for \eqref{EDSKSlin} is very similar with the proof of uniqueness for \eqref{edsKS} which has already been done and we leave it to the reader.
\vskip0.3cm

\textit{Step 5.} It remains to check (ii) to conclude. We thus have to prove the existence of a solution to \eqref{EDSKSlin}. To this aim, we use a Picard iteration. We thus consider the constant process $X_t^0=x_0$ and define recursively
\begin{align*}
X_t^{n+1}=x_0-\chi\int_0^t\int_{\mathbb R^2}K(X_s^n-x)f_s(dx)ds+\sqrt2 B_t.
\end{align*}
Using the same kind of arguments as in the proof of Theorem \ref{exuniedsKS}, we get
\begin{align*}
\mathbb E(\sup_{[0,T]}|X_t^{n+1}-X_t^n|)\leq C\int_0^T\mathbb E[|X_s^n-X_s^{n-1}|](1+||f_s||_{L^p})ds.
\end{align*}
Since $\int_0^T(1+||f_s||_{L^p})ds<\infty$, we classically deduce that $\sum_n \mathbb E(\sup_{[0,T]}|X_t^{n+1}-X_t^n|)<\infty$, so that there is a continuous adapted process $(X_t)_{t \geq0}$ such that for all $T>0$, $\lim_n\mathbb E\big[\sup_{[0,T]}|X_t-X_t^n|\big]=0$. This $L^1$ convergence implies that $(X_t)_{t \geq0}$ is solution to \eqref{EDSKSlin}, which concludes the proof. \hfill$\square$
\vskip0.5cm

The following result ensures that uniqueness holds for \eqref{eq:KS}.
\begin{theo}\label{thm:Wdistance}
Let $f_0$ and $g_0$ be two probability measures with finite first moment. Let $(f_t)_{t\geq0}$ and $(g_t)_{t\geq0}$ be two solutions to \eqref{eq:KS} lying in $L_{loc}^\infty([0,\infty),\textbf{P}_1(\mathbb R^2))\cap L_{loc}^1([0,\infty);L^p(\mathbb R^2))$ for some $p>2/(1-\alpha)$ starting from $f_0$ and $g_0$ respectively. Then
\begin{align*}
\mathcal W_1(f_t,g_t)\leq \mathcal W_1(f_0,g_0)\exp\Big(C\int_0^t(1+||f_s+g_s||_{L^p})ds\Big).
\end{align*}
\end{theo}

\textbf{Proof.}
Let thus $p>2/(1-\alpha)$, $(f_t)_{t\geq0}$ and $(g_t)_{t\geq0}$ be two solutions to \eqref{eq:KS} lying in $L_{loc}^\infty([0,\infty),\textbf{P}_1(\mathbb R^2))\cap L_{loc}^1([0,\infty);L^p(\mathbb R^2))$. For any $s\geq0$, we consider the probability measure $R_s$ on $\mathbb R^2\times\mathbb R^2$ with first (respectively second) marginal equal to $f_s$ (resp. $g_s$) such that 
\begin{align*}
\mathcal W_1(f_s,g_s)=\int_{\mathbb R^2\times\mathbb R^2}|x-y|R_s(dx,dy),
\end{align*}
and we consider $(X_0,Y_0)$ with law $R_0$. We finally set
\begin{align*}
X_t&=X_0-\chi\int_0^t\int_{\mathbb R^2}K(X_s-x)f_s(dx)ds+\sqrt2 B_t,\\
Y_t&=Y_0-\chi\int_0^t\int_{\mathbb R^2}K(Y_s-x)g_s(dx)ds+\sqrt2 B_t.
\end{align*}
Using It\^o's formula, we see that $h$ defined by $h_t:=\mathcal L(X_t)$ satisfies \eqref{KSlinweak} and Lemma \ref{lemuni} ensures us that $\mathcal L(X_t)=f_t$. Similarly, we also have $\mathcal L(Y_t)=g_t$. Using the same arguments as in the proof of Theorem \ref{exuniedsKS}, we easily get
\begin{align*}
\mathbb E(|X_t-Y_t|)\leq \mathbb E[|X_0-Y_0|]+C\int_0^t\mathbb E[
|X_s-Y_s|](1+||f_s+g_s||_{L^p})ds.
\end{align*}
Using the Gr\"onwall's Lemma and recalling that $\mathbb E[|X_0-Y_0|]=\mathcal W_1(f_0,g_0)$, we get
\begin{align*}
\mathbb E(|X_t-Y_t|)\leq \mathcal W_1(f_0,g_0)\exp\Big(C\int_0^t(1+||f_s+g_s||_{L^p})ds\Big),
\end{align*}
which concludes the proof since $\mathcal W_1(f_t,g_t)\leq \mathbb E(|X_t-Y_t|)$. \hfill$\square$
\vskip0.3cm

We can now give the proof of our well-posedness result for \eqref{eq:KS}.

\textbf{Proof of Theorem \ref{thm:ExUnEDP} (i).}
The existence follows by Theorem \ref{exuniedsKS}. Indeed consider $(X_t)_{t\geq0}$ the unique solution of \eqref{edsKS} with initial law $f_0$ and set for $t\geq0$ $f_t:=\mathcal L(X_t)$. Thanks to the Remark \ref{rmq:weaksol}, $f_t$ is a weak solution to \eqref{eq:KS} in the sense given by Definition \ref{def:Weak} and \eqref{Fishmom} is exactly \eqref{eq:FinFirstMom}. 

For uniqueness, consider two weak solutions $(f_t)_{t\geq0}$ and $(g_t)_{t\geq0}$ of \eqref{eq:KS} satisfying \eqref{eq:FinFirstMom} with the same initial condition $f_0\in \textbf{P}_1(\mathbb R^2)$. Then Theorem \ref{thm:Wdistance} ensures that  $\mathcal W_1(f_t,g_t)= 0$ for any $t\geq0$ which concludes the proof. \hfill$\square$
\vskip0.5cm
We end this section with the proof of our propagation of chaos result.

\textbf{Proof of Theorem \ref{propchaos} (i).} We consider $Q^N:=\frac{1}{N}\sum_{i=1}^N\delta_{(X_t^{i,N})_{t\geq0}}$. By Lemma \ref{tension}, the family $\{\mathcal{L}(Q^N),N \geq2\}$ is tight in $\textbf{P}(\textbf{P}(C([0,\infty),\mathbb{R}^2)))$. Furthermore, by proposition \ref{CVversX}, any limit point of $\mathcal{Q}^N$ belongs a.s. to the set of all probability measures $f\in\textbf{P}(C([0,\infty),\mathbb{R}^2)$ such that $f$ is the law of a solution to \eqref{edsKS} satisfying \eqref{eq:Derivative}. But by Theorem \ref{exuniedsKS}, this set is reduced to $\mathcal L((X_t)_{t\geq0})=:f$. We thus deduce that $Q^N$ goes in law to $f$ as $N\rightarrow\infty$ which concludes the proof of $(i)$.

\section{Renormalization and entropic chaos}

In this section, we first deal with the renormalization  which will give us the dissipation of entropy for the solution to \eqref{eq:KS}. From this, we will be able to show the entropic chaos for the system \eqref{systpart}, which will conclude this paper.\\

\textbf{Proof of Theorem \ref{thm:ExUnEDP} (ii).}
We adapt the ideas used in \cite{FHM} for the 2D vortex model to our case, which in particular has a non divergence free kernel. We split the proof in four steps plus a Step 0 which is nothing but direct results of what we have already done. We consider the unique weak solution $f=(f_t)_{t\geq0}$ of \eqref{eq:KS}. In step 1 we deal with the necessary estimates on $K\ast f$ and $\nabla\cdot(K\ast f)$ to regularize $f$. In step 2 we show the convergence of a regular version of $f$ towards $f$. In step 3, we improve the regularity of the solution using a well-known bootstrap argument. Finally, in step 4 we prove the renormalization property.\\

We first observe that by construction, $f$ satisfies \eqref{eq:Fishinff}. Indeed, for any $t\geq0$, we considered $f_t$ as the law of $X_t$, where $(X_t)_{t\geq0}$ is the unique solution to \eqref{edsKS}, obtained by Proposition \ref{CVversX} and Lemma \ref{tension}, so that \eqref{eq:FisherInf} (which englobes~\eqref{eq:Fishinff}) is satisfied. \\

{\it Step 0. Direct Estimates}. We start by noticing that Lemma~\ref{lemLpFish} and~\eqref{eq:Fishinff} implies directly~\eqref{eq:Derivative} and also that for any $p\in[1,\infty)$ and all $T>0$,
\begin{equation}\label{eq:First}
 f\in L^{p/(p-1)}(0,T;L^p(\mathbb R^2)).
\end{equation} 

{\it Step 1. First Estimates}.  The aim of this step is to prove that for any $q>2/\alpha$ and all $T>0$:
\begin{equation}\label{eq:Third}
(K\ast f)\in L^{2q/(\alpha q-2)}(0,T;L^{q}(\mathbb R^2)),
\end{equation}
and
\begin{equation}\label{eq:Fourth}
 \nabla_x\cdot(K\ast f)=K\ast(\nabla_x\cdot f)\in L^{2q/(q(1+\alpha)-2)}(0,T;L^q(\mathbb R^2)).
\end{equation}\vspace{0.1cm}

Let us remember the Hardy-Littlewood-Sobolev inequality in 2D: for $1\leq p<2/(2-\alpha)$,
$$\left\|\,\int_{\mathbb R^2}\frac{f(y)}{|\cdot-y|^{2-(2-\alpha)}}dy\,\right\|_{2p/(2-(2-\alpha)p}\leq C_{\alpha,p}\|f\|_p.$$ Using \eqref{eq:First} we get that for any $p\in(1,2/(2-\alpha))$ and all $T>0$,
\begin{equation}\nonumber
(K\ast f)\in L^{p/(p-1)}(0,T;L^{2p/(2-(2-\alpha)p)}(\mathbb R^2)),
\end{equation}
and under the change of variables $q=2p/(2-(2-\alpha)p)$ we easily deduce \eqref{eq:Third}.

Similarly, but using \eqref{eq:Derivative} instead of \eqref{eq:First}, we get that for any $p\in(1,2/(2-\alpha))$ and all $T>0,$
\begin{equation}\nonumber
 \nabla_x\cdot(K\ast f)\in L^{2p/(3p-2)}(0,T;L^{2p/(2-(2-\alpha)p)}(\mathbb R^2)),
\end{equation}
applying the same change of variables $q=2p/(2-(2-\alpha)p)$ we get \eqref{eq:Fourth}.\\

{\it Step 2. Continuity}. Consider $T>0$ fixed. For $q>2/\alpha$ we have that $2q/(q(1+\alpha)-2)<q/(q-1)$, then using \eqref{eq:First} with $p=q/(q-1)>1$, and \eqref{eq:Fourth}, we get that $f\,\nabla_x\cdot(K\ast f)$ belongs to $L^1(0,T;L^1(\mathbb R^2))$. The following lemma follows directly:
\begin{lemme}\label{DiPerna}
 Consider a mollifier sequence $(\rho_n)$ on $\mathbb R^2$ and introduce the mollified function $f^n_t:= f_t\ast\rho_n$. Clearly, $f_t^n\in C([0,\infty),L^1(\mathbb R^2))$. For all $T>0$, there exists $r^n\in L^1(0,T;L^1_{loc}(\mathbb R^2))$ that goes to $0$ when $n\rightarrow\infty$, and such that
 \begin{equation}\label{eq:Cinfty}
  \partial_t f^n-\chi\nabla_x\cdot((K\ast f)f^n)-\triangle_x f^n=r^n.
\end{equation}
\end{lemme}

\begin{rmq}
The proof of the previous lemma is a modification of \cite[Lemma II.1.(ii) and Remark 4]{dPL}. In fact, for all $T>0$, $f\in L^\infty(0,T;L^1(\mathbb R^2))$ and for any $p>2/\alpha$, $(K\ast f)\in L^1(0,T;L^p(\mathbb R^2))$. That suffices for the existence of $r^n$ given by
\begin{align*}
 &r^n:=\chi\big[\big(\nabla\cdot((K\ast f)f)\big)\ast \rho^n-\nabla\cdot\big((K\ast f)f^n\big)\big],
\end{align*}
which goes to 0 if $n\rightarrow\infty$ in $L^1(0,T;L^1_{loc}(\mathbb R^2))$.
\end
{rmq}

As a consequence of Lemma \ref{DiPerna}, the chain rule applied to the smooth $f^n$ reads
\begin{align}  \label{eq:Beta}
\partial_t\beta(f^n)=&\chi\,\left[(K\ast f)\cdot\nabla_x\beta(f^n)+\beta^\prime(f^n)f^n\nabla_x\cdot(K\ast f)\right]\\
 &+\triangle_x\beta(f^n)-\beta^{\prime\prime}(f^n)|
\nabla_xf^n|^2+\beta^\prime(f^n)r^n,\nonumber
\end{align}
for any $\beta\in C^1(\mathbb R)\cap W_{loc}^{2,\infty}(\mathbb R)$ such that $\beta^{\prime\prime}$ is piecewise continuous and vanishes outside of a compact set. Since the equation \eqref{eq:Cinfty} with $(K\ast f)$ fixed is linear in $f^n$, the difference $f^{n,k}:=f^n-f^k$ satisfies \eqref{eq:Cinfty} with $r^n$ replaced by $r^{n,k}:=r^n-r^k\rightarrow0$ in $L^1(0,T;L_{
loc}^1(\mathbb R^2))$ and then also \eqref{eq:Beta} (with again $f^n$ and $r^n$ changed in $f^{n,k}$ and $r^{n,k}$).

Now, choosing $\beta(s)=\beta_1(s)$ where $\beta_1(s)=s^2/2$ for $|s|\leq1$ and $\beta_1(s)=|s|-1/2$ for $|s|\geq 1$. It is clear that $\beta\in C^1(\mathbb R)$, that $\beta^\prime,\beta^{\prime\prime}\in L^\infty(\mathbb R)$ and that the 
second derivative has compact support. For any non-negative $\psi\in C_c^2(\mathbb R^2)$, we obtain
\begin{align*}
& \frac{d}{dt}\int_{\mathbb R^2}\beta_1(f^{n,k}(t,x))\psi(x)\,dx\\
&\quad =\int_{\mathbb R^2}\chi\,\left[(K\ast f)\cdot\nabla_x\beta_1(f^{n,k})+\beta_1^\prime(f^{n,k})f^{n,k}\nabla_x\cdot(K\ast f)\right]\psi(x)\,dx\\
& \quad\quad+\int_{\mathbb R^2}\left[\triangle_x\beta_1(f^{n,k})-\beta_1^{\prime\prime}(f^{n,k})|\nabla_xf^{n,k}|^2+\beta_1^\prime(f^{n,k})r^{n,k}\right]\psi(x)\,dx\\
& \quad\leq\int_{\mathbb R^2}\left|r^{n,k}(t,x)\right|\psi(x)\,dx+\int_{\mathbb R^2}\beta_1(f^{n,k})\triangle_x\psi\,dx\\
&\quad\quad+\chi\int_{\mathbb R^2}|f^{n,k}\,\nabla_x\cdot(K\ast f)|\psi(x)\,dx-\chi\int_{\mathbb R^2}\beta_1(f^{n,k})\nabla_x\cdot\big((K\ast f)\psi(x)\big)\,dx,
\end{align*}
where we have used that $|\beta_1^\prime|\leq1$ and that $\beta_1^{\prime\prime}\geq0$. We know that $f_0\in L^1(\mathbb R^2)$ then $f^{n,k}(0)\rightarrow0$ in $L^1(\mathbb R^2)$, also that $r^{n,k}\rightarrow0$ in $L^1(0,T;L^1_{loc}(\mathbb R^2))$. It is not difficult to see that $\beta_1(f^{n,k})(K\ast f)\rightarrow0$ in $L^1(0,T;L^1_{loc}(\mathbb R^2))$, (because $\beta_1$ is sub-linear, and for all $0<\alpha<1$ there is $q:=p/(p-1)>2/\alpha$, then using \eqref{eq:First} and \eqref{eq:Third}: $f^{n,k}\rightarrow0
\text{ in }L^{p/(p-1)}(0,T;L^p(\mathbb R^2)),$ and $(K\ast f)\in L^{q/(q-1)}(0,T;L^q(\mathbb R^2))).$

The same arguments apply to $\beta_1(f^{n,k})\nabla_x\cdot(K\ast f)$ and $|f^{n,k}\,\nabla_x\cdot(K\ast f)|$, and then both goes to $0$ as $n,k\rightarrow\infty$ in $L^1(0,T;L^1_{loc}(\mathbb R^2))$. Finally, we get
\begin{equation}\nonumber
 \sup_{t\in[0,T]}\int_{\mathbb
 R^2}\beta_1(f^{n,k}(t,x))\psi(x)\,dx\xrightarrow[n,k\rightarrow\infty]{}0.
\end{equation}

Since $\psi$ is arbitrary, we deduce that there exists $\bar f\in C([0,\infty);L^1_{loc}(\mathbb R^2))$ so that $f^n\rightarrow\bar f$ in $C([0,\infty);L^1_{loc}(\mathbb R^2))$ with the topology of the uniform convergence on any compact subset in time. Together with the convergence $f^n\rightarrow f$ in $C([0,\infty);{\bf P}(\mathbb R^2))$ we get that $f=\bar f$. We end this Step by concluding that, with the same convention for the notion of convergence on $[0,\infty)$: $f^n\rightarrow f$ in $C([0,\infty);L^1(\mathbb R^2)).$\\

{\it Step 3. Additional estimates}. From \eqref{eq:Beta}, we know that for all $0<t_0<t_1$, all $\psi\in C_c^2(\mathbb R^2)$,
\begin{align}
\label{eq:Beta2} &\int_{\mathbb R^2}\beta(f^n_{t_1})\psi(x)\,dx+\int_{t_0}^{t_1}\int_{\mathbb R^2}\beta^{\prime\prime}(f_s^n)|\nabla_x f_s^n|^2\psi(x)\,dx\,ds\\
\nonumber &\quad=\int_{\mathbb R^2}\beta(f^n_{t_0})\psi(x)\,dx+\int_{t_0}^{t_1}\int_{\mathbb R^2}\beta^\prime(f^n_s)r^n\psi(x)\,dx\,ds\\
\nonumber &\qquad+\int_{t_0}^
{t_1}\int_{\mathbb R^2}\beta(f^n_s)\big[\triangle_x\psi(x)-\chi\,(K\ast f)\nabla_x\psi(x)\big]\,dx\,ds\\
\nonumber &\qquad+\chi\int_{t_0}^{t_1}\int_{\mathbb R^2}\big[\beta^\prime(f^{n}_s)f^{n}_s-\beta(f^n_s)\big]\psi(x)\,\nabla_x\cdot(K\ast f)\,dx\,ds.
\end{align}
Let us choose $0\leq\psi\in C_c^2(\mathbb R^2)$ and $\beta\in C^1(\mathbb R)\cap W^{2,\infty}_{loc}(\mathbb R)$ convex such that $\beta^{\prime\prime}$ is non-negative and vanishes outside of a compact set (notice that, there is a constant $C>0$ such that $s\beta^\prime(s)\leq C\beta(s)$).
We can pass to the limit as $n\rightarrow\infty$ (for details see Step 2) to get
\begin{align*}
\int_{\mathbb R^2}\beta(f_{t_1})\psi(x)\,dx\leq&\int_{\mathbb R^2}\beta(f_{t_0})\psi(x)\,dx\\
&+\int_{t_0}^{t_1}\int_{\mathbb R^2}\beta(f_s)\left[\triangle_x\psi(x)-\chi(K\ast f)\nabla_x\psi(x)\right]\,dx\,ds\\
&+\chi\int_{t_0}^{t_1}\int_{\mathbb R^2}\left[-\beta(f_s)+\beta^\prime(f_s)f_s\right]\psi(x)\,\nabla_x\cdot(K\ast f)\,dx\,ds.
\end{align*}

It is not hard to deduce, by approximating $\psi\equiv1$ by a well-chosen sequence $\psi_R$ that
\begin{equation}\nonumber
\nonumber \int_{\mathbb R^2}\beta(f_{t_1})\,dx\leq\int_{\mathbb R^2}\beta(f_{t_0})\,dx+\chi\int_{t_0}^{t_1}\int_{\mathbb R^2}\left[-\beta(f_s)+\beta^\prime(f_s)f_s\right]\,\nabla_x\cdot(K\ast f)\,dx\,ds.
\end{equation}
whenever $\beta$ is admissible.\\

Now we deal with the regularity in space of \eqref{eq:Regularity}. Let us start by noticing that taking $p>2/(1-\alpha)$:
\begin{equation}
\nabla_x(K\ast f)(x)=\int_{\mathbb R^2}\frac{(1-\alpha)f(y)}{|x-y|^{1+\alpha}}\,dy,\label{eq:NabInf}
\end{equation}
so that using \eqref{Jalpha},
\begin{equation}\nonumber
 \int_0^T\|\nabla_x(K\ast f_s)\|_{L^\infty(\mathbb R^2)}\leq C(\alpha,p)\int_0^T\left(\|f_s\|_{L^p(\mathbb R^2)}+1\right)<\infty,
\end{equation}
and due to the fact that $s\beta^\prime(s)\leq C\beta(s)$, we get
\begin{eqnarray*}
\nonumber \int_{\mathbb R^2}\beta(f_{t_1})\,dx&\leq&
\int_{\mathbb R^2}\beta(f^n_{t_0})\,dx\\
&&+(C+1)\chi\int_{t_0}^{t_1}\|\nabla_x(K\ast f)(x)\|_{L^\infty(\mathbb R^2)}\int_{\mathbb R^2}\beta(f_s)\,dx\,ds.
\end{eqnarray*}
Then Gr\"onwall's lemma implies that for all $0<t_0<t_1<T$,
\begin{eqnarray*}
\nonumber \int_{\mathbb R^2}\beta(f_{t_1})\,dx&\leq&C(\alpha,T)\int_{\mathbb R^2}\beta(f^n_{t_0})\,dx.
\end{eqnarray*}
Finally letting $\beta(s)\rightarrow|s|^q/q$, we get that for all $q\geq1$ and all $0<t_0<t_1<T$,
\begin{equation}\label{eq:q}
 \|f(t_1,\cdot)\|_{L^q(\mathbb R^2)}\leq C(q,\alpha,T)\|f(t_0,\cdot)\|_{L^q(\mathbb R^2)}.
\end{equation}

Coming back to \eqref{eq:Beta2} and using $\beta_M(s)=s^2/2$ for $|s|\leq1$ and $\beta_M(s)=M|s|-M^2/2$ for $|s|\geq M$, we have
\begin{eqnarray*}&&\int_{\mathbb R^2}\beta_M(f^n_{t_1})\psi\,dx+\int_{t_0}^{t_1}\int_{\mathbb R^2}\one_{|f_s|\leq M}|\nabla_x f_s^n|^2\psi\,dx\,ds\\
&&=\int_{\mathbb R^2}\beta_M(f^n_{t_0})\psi\,dx+\int_{t_0}^{t_1}\int_{\mathbb R^2}\beta_M^\prime(f^n_s)r^n\psi(x)\,dx\,ds\\
&&\quad\int_{t_0}^{t_1}\int_{\mathbb R^2}\beta_M(f^n_s)(K\ast w)\big[\triangle\psi(x)-\,\chi\nabla_x\psi(x)\big]\,dx\,ds\\
&&\quad+\chi\int_{t_0}^{t_1}\int_{\mathbb R^2}\big[\beta_M^\prime(f^{n}_s)f^{n}_s-\beta_M(f^n_s)\big]\psi(x
)\,\nabla_x\cdot(K\ast f)\,dx\,ds,
\end{eqnarray*}
similarly as above we first make $n\rightarrow\infty$, then we approximate $\psi\equiv1$ by a well-chosen sequence $\psi_R$ and make $R\rightarrow\infty$, and finally make the limit $M\rightarrow\infty$ to find that for 
every $T\geq t_1\geq t_0\geq0$:
\begin{align*}
&\int_{\mathbb R^2}|f_{t_1}|^2\,dx+\int_{t_0}^{t_1}\int_{\mathbb R^2}|\nabla_x f_s|^2\,dx\,ds\\
&\quad\leq\int_{\mathbb R^2}|f_{t_0}|^2\,dx+\chi\int_{t_0}^{t_1}\|\nabla_x(K\ast f)(x)\|_{L^\infty(\mathbb R^2)}\int_{\mathbb R^2}|f_s|^2\,dx\,ds.
\end{align*}
We conclude, using \eqref{eq:q}, that for all $0<t_0<T$ and any $q\in[1,\infty)$:
\begin{equation}\label{eq:wInf}
 f\in L^\infty(t_0,T;L^q(\mathbb R^2))\quad\text{and}\quad\nabla_xf\in L^2((t_0,T)\times\mathbb R^2).
\end{equation}\vspace{0.1cm}

To get the continuity in time of \eqref{eq:Regularity}, we need to improve even more the estimates on $f$ which will be achieved using a bootstrap argument. First, fixing $p>2/(2-\alpha)$ we notice that for all $t_0>0$ $$\|K\ast f_t\|_{L^\infty}\leq C(p)(1+\|f_t\|_{L^p})\Rightarrow K\ast f_t\in L^\infty(t_0,T;L^\infty(\mathbb R^2)),$$ and thanks to \eqref{eq:NabInf} and \eqref{eq:wInf}: $$\|\nabla_x(K\ast f_t)\|_{L^\infty}\leq C(p)(1+\|f_t\|_{L^p})\Rightarrow \nabla_x(K\ast f_t)\in L^\infty(t_0,T;L^\infty(\mathbb R^2)),$$ we thus have
\begin{equation}\nonumber
 \partial_t f-\triangle_x f =\big[\chi f\,\nabla_x\cdot (K\ast f)+(K\ast f)\cdot\nabla_x f\big] \in L^2((t_0,T)\times\mathbb R^2),
\end{equation}
and \cite[Theorem X.11]{Brezis} provides the maximal regularity in $L^2$ spaces for the heat equation, in other words: for all $t_0>0$
\begin{equation}\nonumber
f\in L^\infty(t_0,T;H^1(\mathbb R^2))\cap L^2(t_0,T;H^2(\mathbb R^2).
\end{equation}

\begin{rmq}
We emphasize that the previous bound is true for all $t_0$. In fact, when $f_{t_0}\in H^1(\mathbb R^2)$, the maximal regularity implies the above bound in the time interval $[t_0,\infty)$. But thanks to~\eqref{eq:wInf}, we can find $t_0$ arbitrary close to $0$ such that $f_{t_0/2}\in H^1(\mathbb R^2)$, then we get the conclusion.
\end{rmq}

Using now the interpolation inequality, there exists a constant $C>0$ such that
\begin{equation}\nonumber
 \|\nabla_x f\|_{L^3(\mathbb R^2)}\leq C\|D^2 f\|_{L^2(\mathbb R^2)}^{2/3}\|f\|_{L^2(\mathbb R^2)}^{1/3},
\end{equation}
which implies
\begin{equation}\nonumber
 \int_{t_0}^T\|\nabla_x f\|_{L^3(\mathbb R^2)}^3\,ds\leq
 C\int_{t_0}^T\|D^2 f\|_{L^2(\mathbb R^2)}^{2}\|f\|_{L^2(\mathbb R^2)}<\infty.
\end{equation}

Thanks to the previous calculus and again \cite[Theorem X.12]{Brezis} we conclude that $\partial_t f,\nabla_x f\in L^3((t_0,T)\times\mathbb R^2)$ and then the Morrey's inequality implies that for all $t_0>0$
\begin{equation}\nonumber
 f\in C^0((t_0,T)\times\mathbb R^2),
\end{equation}
all together allow us to deduce that
\begin{equation}\nonumber
 f\in C([0,T);L^1(\mathbb R^2))\cap C((0,T);L^2(\mathbb R^2)).
\end{equation}

We can go even further iterating this argument, using the interpolation inequality and the Sobolev inequality, to deduce that $\nabla_x f\in L^p((t_0,T)\times\mathbb R^2)$ for any $1< p<\infty$,  $\left[\chi f\,\nabla_x\cdot (K\ast f)+(K\ast f)\cdot\nabla_x f\right] \in L^p((t_0,T)\times\mathbb R^2)$ for all $t_0>0$. Then the maximal regularity of the heat equation in $L^p$ spaces (see \cite[Theorem X.12]{Brezis}) implies that for all $t_0>0$
\begin{equation}\nonumber
 \partial_tf,\nabla_x f\in L^p((t_0,T)\times\mathbb R^2),
\end{equation}
and then using again the Morrey's inequality: $f\in C^{0,\alpha}((t_0,T)\times\mathbb R^2)$ for any $0<\alpha<1$, and any $t_0>0$. All together allow us to conclude \eqref{eq:Regularity}.\vspace{0.3cm}

{\it Step 4. Renormalization}. To end the proof we show \eqref{eq:Renormalization}. Let thus $\beta\in C^1(\mathbb R)\cap W^{2,\infty}_{loc}(\mathbb R)$ sub-linear, such that $\beta^{\prime\prime}$ is piecewise continuous and vanishes outside of a compact set. Thanks to \eqref{eq:wInf}, we can pass to the limit in the similar identity as \eqref{eq:Beta2} obtained for time dependent test functions $\psi\in C_c^2([0,\infty)\times\mathbb R^2)$ to get
\begin{align}\label{eq:Beta3}
&\int_{t_0}^\infty\int_{\mathbb R^2}\beta^{\prime\prime}(f_s)|\nabla_x f_s|^2\psi_s\,dx\,ds=\int_{\mathbb R^2}\beta(f_{t_0})\psi_{t_0}\,dx\\
\nonumber &\quad-\chi\int_{t_0}^\infty\int_{\mathbb R^2}\psi_s(x)\,\nabla_x\cdot(K\ast f)\big(f_s\beta^\prime(f_s)-\beta(f_s)\big)\,dx\,ds\\
\nonumber &\qquad\qquad+\int_{t_0}^\infty\int_{\mathbb R^2}\beta(f_s)\big(\triangle_x\psi_s(x)-(K\ast f)\nabla_x\psi_s(x)+\partial_t\psi_s(x)\big)\,dx\,ds.
\end{align}
In the case $\psi\geq0$ and $\beta^{\prime\prime}\geq0$ we can pass to the limit $t_0\rightarrow0$ thanks to monotonous convergence in the first term,  the continuity property obtained in Step 2 in the second term, and the monotonous convergence in the other terms (recall that $s\beta^\prime(s)\leq\beta(s)$, $\beta$ is sub-linear and $|f|(1+|K\ast f|+|\nabla\cdot(K\ast f)|)$ belongs to $L^1(0,T;L^1(\mathbb R^2)$ thanks to \eqref{eq:Third} and \eqref{eq:Fourth}). We get 
\begin{align}\label{eq:Beta4}
&\int_{0}^\infty\int_{\mathbb R^2}\beta^{\prime\prime}(f_s)|\nabla_x f_s|^2\psi_s\,dx\,ds=\int_{\mathbb R^2}\beta(f_{0})\psi_{t_0}\,dx\\
\nonumber &\qquad+\int_{0}^\infty\int_{\mathbb R^2}\beta(f_s)\left[\triangle_x\psi_s-\chi\nabla_x((K\ast f)\cdot\psi_s)+\partial_t\psi_s\right]\,dx\,ds\\
\nonumber &\qquad\qquad+\chi\int_{0}^\infty\int_{\mathbb R^2}\beta^\prime(f_s)f_s\psi_s(x)\,\nabla_x\cdot(K\ast f)\,dx\,ds,
\end{align}
and the bound given by \eqref{eq:Beta4} implies directly that we can pass to the limit $t_0\rightarrow0$ in the general case for $\psi$ in \eqref{eq:Beta3} which is nothing but \eqref{eq:Renormalization} in the distributional sense. \hfill $\square$
\vskip0.5cm

We now give a useful lemma for the entropic chaos.

\begin{lemme} \label{lementrofish}
Let $\alpha\in(0,1)$ and $f_0\in \textbf{P}_1(\mathbb R^2)$ such that $H(f_0)<\infty$. Let $(f_t)_{t\geq0}$ be the unique solution of \eqref{eq:KS} satisfying \eqref{eq:FinFirstMom}. Then
\begin{equation}\label{eq:FuncEq}
H(f_t)+\int_0^t I(f_s)ds=H(f_0)+\chi(1-\alpha)\int_0^t\int_{\mathbb R^2}\int_{\mathbb R^2} \frac{f_s(dx)f_s(dy)}{|x-y|^{\alpha+1}}ds.
\end{equation}
\end{lemme}

\textbf{Proof.} For $m>1$, let us take $\beta_m\in C^1(\mathbb R)\cap W^{2,\infty}_{loc}(\mathbb R)$ given by
\begin{equation}\nonumber
\beta_m(s)=\begin{cases}
 s\log(s)+(1-s)/m & \text{for }m^{-1}\leq s\leq m,\\
 \beta_m(m_-)+\beta_m^\prime(m_-)(s-m) & \text{for }s>m,\\
 \beta_m(m^{-1}_+)+\beta_m^\prime(m^{-1}_+)\big(s-m^{-1})\big) & \text{for }s<m^{-1},
\end{cases}
\end{equation}
so that $\beta_m(s)\leq Cs$ and $\beta_m\rightarrow s\log(s)$ for any $s>0$.

Since $\beta_m$ is admissible (in the sense of Theorem~\ref{thm:ExUnEDP}), then using \eqref{eq:Renormalization} we get that for any $\psi\in C_c^\infty(\mathbb R^2)$,  
\begin{align*}
  \int \beta_m(f_t)\psi\,dx-\int\beta_m(f_{0})\psi\,dx=&\chi\,\int_{0}^t\int\nabla_x\cdot(K\ast f)\big(f\beta_m^\prime(f)-\beta_m(f)\big)\psi\,dx\,ds\\
  &+\int_{0}^t\int\beta_m(f)\big(\triangle_x\psi-\chi(K\ast f)\cdot\nabla_x\psi\big)\,dx\,ds\\
  &-\int_{0}^t\int\beta_m^{\prime\prime}(f)|\nabla_x f|^2\psi\,dx\,ds,
\end{align*}
using that $\beta_m^{\prime\prime}(s)$ is non-negative, that $\beta_m$ growths linearly at $+\infty$ and that $(f_s)_{s\geq0}$ is non-negative we can make $\psi\rightarrow 1$ to get
\begin{align*}
  \int \beta_m(f_t)\,dx-\int\beta_m(f_{0})\,dx=&\,\chi\,\int_{0}^t\int\nabla_x\cdot(K\ast f)\big(f\beta_m^\prime(f)-\beta_m(f)\big)\,dx\,ds\\
  &-\int_{0}^t\int\beta_m^{\prime\prime}(f)|\nabla_x f|^2\,dx\,ds.
\end{align*}
In fact, the first and the second terms converge thanks to monotonous convergence and that $|\beta_m(s)|\leq C|s|$. The third term is a consequence of the monotonous convergence, that $\beta_m^\prime(s)$ is bounded, and that $f\,\nabla\cdot(K\ast f)$ (resp. $|f(K\ast f)|$ for the fourth term) is integrable by \eqref{eq:Fourth} (resp. \eqref{eq:Third}). The last term is a consequence of \eqref{eq:FisherInf}.\\

Finally, we notice that in the interval $(0,1]$ the function $-\beta_m$ increases to $-s\log(s)$ while in the interval $[1,\infty)$, $\beta_m(s)$ increases to $s\log(s)$. Thanks to the monotonous convergence we can make $m\rightarrow\infty$ and using the integrability of all the limits we get \eqref{eq:FuncEq}.
\hfill $\square$

\vskip0.5cm

It remains to conclude with the proof of the entropic chaos. 

\textbf{Proof of Theorem \ref{propchaos} (ii).}
We only have to prove that for each $t\geq0$, $H(F_t^N)$ tends to $H(f_t)$. To this aim, we first show that for 
any $t\geq0$
\begin{align} \label{limsupH+I}
L:=\limsup_N\Big[H(F_t^N)+\int_0^tI(F_s^N)ds\Big]\leq H(f_t)+\int_0^tI(f_s)ds.
\end{align}
Let $t\geq0$ be fixed. Using \eqref{entroFishFt} and recalling that $H(F_0^N)\rightarrow H(f_0)$ by assumption, we have
\begin{align*}
L\leq H(f_0)+\limsup_N \frac{\chi(1-\alpha)}{N^2}\sum_{i\neq j}\int_0^t\mathbb E\Big[\frac{1}{|X_s^{i,N}-X_s^{j,N}|^{\alpha+1}}\Big]ds,
\end{align*}
so that using that $H(f_t)+\int_0^tI(f_s)ds=H(f_0)+\chi(1-\alpha)\int_0^t\int_{\mathbb R^2}\int_{\mathbb R^2}\frac{f_s(dx)f_s(dy)}{|x-y|^{\alpha+1}}ds$ by Lemma \ref{lementrofish}, we only have to prove that
\begin{align*}
 \lim_{N\rightarrow\infty}\frac{1}{N^2}\int_0^t\mathbb E\Big[\sum_{i\neq j}\frac{1}{|X_s^{i,N}-X_s^{j,N}|^{\alpha+1}}\Big]ds=\int_0^t\int_{\mathbb R^2}\int_{\mathbb R^2}\frac{f_s(dx)f_s(dy)}{|x-y|^{\alpha+1}}ds.
\end{align*}
By exchangeability, it suffices to prove that, as $N\rightarrow\infty$,
\begin{align*}
D_N:=\int_0^t\mathbb E\Big[\frac{1}{|X_s^{1,N}-X_s^{2,N}|^{\alpha+1}}\Big]ds\rightarrow\int_0^t\int_{\mathbb R^2}\int_{\mathbb R^2}\frac{f_s(dx)f_s(dy)}{|x-y|^{\alpha+1}}ds=:D.
\end{align*}
For any $\epsilon>0$, we have
\begin{align*}
|D-D_N|\leq |D-D_\epsilon|+|D_\epsilon-D_{N,\epsilon}|+|D_{N,\epsilon}-D_N|,
\end{align*}
where $D_{N,\epsilon}=\int_0^t\mathbb E\Big[\frac{1}{(|X_s^{1,N}-X_s^{2,N}|\vee\epsilon)^{\alpha+1}}\Big]ds$ and $D_\epsilon=\int_0^t\int_{\mathbb R^2}\int_{\mathbb R^2}\frac{f_s(dx)f_s(dy)}{(|x-y|\vee\epsilon)^{\alpha+1}}ds$. Using that for any $\epsilon>0$ fixed, the function $(x,y)\mapsto (|x-y|\vee\epsilon)^{-\alpha-1}$ is bounded continuous and that $\mathcal L(X_s^{1,N},X_s^{2,N})$ goes weakly to $f_s\otimes f_s$ 
for any $s\geq0$, we have $\lim_N \mathbb E\Big[\frac{1}{(|X_s^{1,N}-X_s^{2,N}|\vee\epsilon)^{\alpha+1}}\Big]=\int_{\mathbb R^2}\int_{\mathbb R^2}\frac{f_s(dx)f_s(dy)}{(|x-y|\vee\epsilon)^{\alpha+1}}$. By dominated convergence, we thus get that $\lim_N|D_\epsilon-D_{N,\epsilon}|=0$. We thus have
\begin{align*}
\limsup_N|D-D_N|\leq |D-D_\epsilon|+\limsup_N|D_{N,\epsilon}-D_N|\quad\forall\epsilon>0.
\end{align*}
Let $\tilde \alpha$ be such that $\alpha+1<\tilde \alpha<2$. We have
\begin{align*}
|D-D_\epsilon|&\leq 2\int_0^t\int_{\mathbb R^2}\int_{\mathbb R^2}\frac{f_s(dx)f_s(dy)}{|x-y|^{\alpha+1}}\one_{\{|x-y|<\epsilon\}}ds\\
&\leq  2\epsilon^{\tilde \alpha-\alpha-1}\int_0^t\int_{\mathbb R^2}\int_{\mathbb R^2}\frac{f_s(dx)f_s(dy)}{|x-y|^{\tilde\alpha}}ds\\
&\leq C\epsilon^{\tilde \alpha-\alpha-1}\int_0^t(1+I(f_s))ds\leq C(1+t)\epsilon^{\tilde \alpha-\alpha-1},
\end{align*}
by Lemma \ref{lemFish} (applied with $F=f_s\otimes f_s$, for which $I(F_s)=I(f_s)$) and \eqref{eq:Fishinff}. Using the same arguments, we also have for any $N\geq2$,
\begin{align*}
|D_{N,\epsilon}-D_N|\leq C\epsilon^{\tilde \alpha-\alpha-1}\int_0^t(1+I(F_s^N))ds\leq C(1+t)\epsilon^{\tilde \alpha-\alpha-1}.
\end{align*}
We thus get that $\limsup_N |D-D_N|=0$ and \eqref{limsupH+I} is proved.

Using \cite[Theorem 3.4 and Theorem 5.7]{HM}, we have
\begin{align} \label{liminfH,I}
\liminf_N H(F_t^N)\geq H(f_t) \quad{\rm and}\quad \liminf_N \int_0^tI(F_s^N)ds\geq \int_0^tI(f_s)ds.
\end{align}
Using \eqref{limsupH+I} and \eqref{liminfH,I}, we easily conclude that 
\begin{align*}
\lim_N H(F_t^N)= H(f_t) \quad{\rm and}\quad \lim_N \int_0^tI(F_s^N)ds= \int_0^tI(f_s)ds,
\end{align*}
which concludes the proof. \hfill$\square$


\begin{thebibliography}{99}
\bibitem{BK}
{Bhatt, A; Karandikar, R. Invariant measures and evolution equations for Markov processes characterized via martingale problems. Ann. Probab.  21  (1993),  no. 4, 2246-2268.}
\bibitem{Brezis}
{Brezis, H. Analyse fonctionnelle. Collection Math\'ematiques Appliqu\'ees pour la Ma\^itrise. Masson, Paris, 1983. Th\'eorie et applications.}
\bibitem{BRE}
{Brezis, H. Convergence in D′ and in L1 under strict convexity. Boundary value problems for partial differential equations and applications, 43-52, RMA Res. Notes Appl. Math., 29, Masson, Paris, 1993.}
\bibitem{BDP}
{Blanchet, A.; Dolbeault, J.; Perthame, B. Two-dimensional Keller-Segel model: optimal critical mass and qualitative properties of the solutions. Electron. J. Differential Equations 2006, No. 44, 32 pp.} 
\bibitem{Calvez13}
{Calvez, V.; Corrias, L. Blow-up dynamics of self-attracting diffusive particles driven by competing convexities. arXiv:1301.7075.}
\bibitem{CCLLV}
{Carlen, E.; Carvalho, M.; Le Roux, J.; Loss, M.; Villani, C. Entropy and chaos in the Kac model. Kinet. Relat. Models  3  (2010),  no. 1, 85-122.}
\bibitem{dPL}
{DiPerna, R., Lions, P.-L. Ordinary differential equations, transport theory and Sobolev spaces, Invent. Math. 98, 3 (1989), 511-547.}
\bibitem{DUR}
{Durrett, R. Stochastic calculus. A practical introduction. Probability and Stochastics Series. CRC Press, Boca Raton, FL, 1996. x+341 pp.}
\bibitem{FHM}
{Fournier, N.; Hauray, M.; Mischler,  S. Propagation of chaos for the 2D viscous vortex model. arXiv:1212.1437} 
\bibitem{Schmeiser1}
{Haskovec, J.; Schmeiser, C. Stochastic Particle Approximation to the Global Measure Valued Solutions of the Keller-Segel model in 2D, J. Stat. Phys. 135 (2009), pp. 133-151.}
\bibitem{Schmeiser2}
{Haskovec, J.; Schmeiser, C. Convergence analysis of a stochastic particle approximation for measure valued solutions of the 2D Keller-Segel system, Comm. PDE 36 (2011), pp. 940-960.}
\bibitem{HM}
{Hauray, M.; Mischler, S. On Kac's chaos and related problems. arXiv:1205.4518.}
\bibitem{HOR1}
{Horstmann, D. From 1970 until present: the Keller-Segel model in chemotaxis and its consequences. I. Jahresber. Deutsch. Math.-Verein. 105 (2003), no. 3, 103-165.}
\bibitem{HOR2}
{Horstmann, D. From 1970 until present: the Keller-Segel model in chemotaxis and its consequences. II. Jahresber. Deutsch. Math.-Verein. 106 (2004), no. 2, 51-69.}
\bibitem{KAC}
{Kac, M. Foundations of kinetic theory.  Proceedings of the Third Berkeley Symposium on Mathematical Statistics and Probability, 1954-1955, vol. III,  pp. 171-197. University of California Press, Berkeley and Los Angeles, 1956.}
\bibitem{KS1}
{Keller, E.F.; Segel,L.A. Initiation of slime mold aggregation viewed as an instability, J. Theoret. Biol., 26 (1970), pp. 399-415.}
\bibitem{KS2}
{Keller, E.F.; Segel, L.A. A model for chemotaxis, J. Theoret. Biol., 30 (1971), pp. 225-234.}
\bibitem{Poupaud}
{Poupaud, F. Diagonal defect measures, adhesion dynamics and Euler equations, Meth. Appl. Anal. 9 (2002), 20pp. 533-561.}
\bibitem{Stevens1}
{Stevens, A. A stochastic cellular automaton, modeling gliding and aggregation of myxobacteria, SIAM J. Appl. Math. Vol. 61 (2000), pp. 172-182.}
\bibitem{Stevens2}
{Stevens, A. The derivation of chemotaxis equations as limit dynamics of moderately interacting stochastic many-particle systems. SIAM J. Appl. Math. Vol. 61 (2000), pp. 183-212.}
\bibitem{SZN}
{Sznitman, A.S. Topics in propagation of chaos. \'Ecole d'\'Et\'e de Probabilit\'es de Saint-Flour XIX—1989, 165–251, Lecture Notes in Math., 1464, Springer, Berlin, 1991.}
\bibitem{TAK}
{Takanobu, S. On the existence and uniqueness of SDE describing an $n$-particle system interacting via a singular potential. Proc. Japan Acad. Ser. A Math. Sci.  61  (1985), no. 9, 287-290.}
\bibitem{VIL2} 
{Villani, C.: Topics in optimal transportation. Graduate Studies in 
Mathematics, 58. American Mathematical Society, Providence, RI, 2003}.
\end{thebibliography}
\end{document}